\documentclass{article}
\usepackage{amsfonts}
\usepackage{amsmath}
\usepackage{graphics}
\usepackage{latexsym}
\usepackage[usenames]{color}

\usepackage{enumitem}
\usepackage{epsfig}

\usepackage{amssymb}

\usepackage{amscd}
\usepackage{psfrag}

\usepackage{amsthm}
\graphicspath{{./plots_for_papers/}}

\def\C{\mathbb{C}}

\def\m{\mathcal{M}_{\Gamma}}
\def\F{\mathcal{F}}

\def\A{\mathcal{A}}
\def\D{\mathbb{D}}

\def\R{\mathbb{R}}

\def\Z{\mathbb{Z}}

\def\H{\mathbb{H}}

\begin{document}
\newtheorem{thm}{Theorem}

\newtheorem{claim}{Claim}

\newtheorem{lemma}{Lemma}
\newtheorem{defn}{Definition}

\newtheorem{cor}{Corollary}

\newtheorem{prop}{Proposition}

\newcommand{\bx}[1]{\bibitem{#1}}

\newtheorem{example}{Example}
\newtheorem{definition}{Definition}
\newtheorem{conj}{Conjecture}
\newtheorem{remark}{Remark}

\def\F{\mathcal{F}}
\title{Dynamics of modular matings}

\author{Shaun Bullett and Luna Lomonaco}

\maketitle

\begin{abstract}
We develop a dynamical theory for the family of holomorphic correspondences $\F_a$ proved by the current authors
to be matings between the modular group and parabolic rational maps in the Milnor slice $Per_1(1)$ (\cite{BL1}). 
Such a mating endows the complement of the limit set of $\F_a$ with the geometry of the hyperbolic plane, equipped with the action of the modular group.  
We introduce bi-infinite coding sequences for geodesics in this complement, utilising continued fraction expressions of end points; we prove landing theorems
for periodic and preperiodic geodesics, and we establish a stronger Yoccoz inequality 
for repelling fixed points of these correspondences than Yoccoz's classical inequality for quadratic polynomials. 
We deduce that the connectedness locus of the family $\F_a$ is contained in a particular lune in parameter space.
\end{abstract}

\section{Introduction}\label{intro}
\begin{figure}
\begin{center}
\scalebox{.29}{\includegraphics{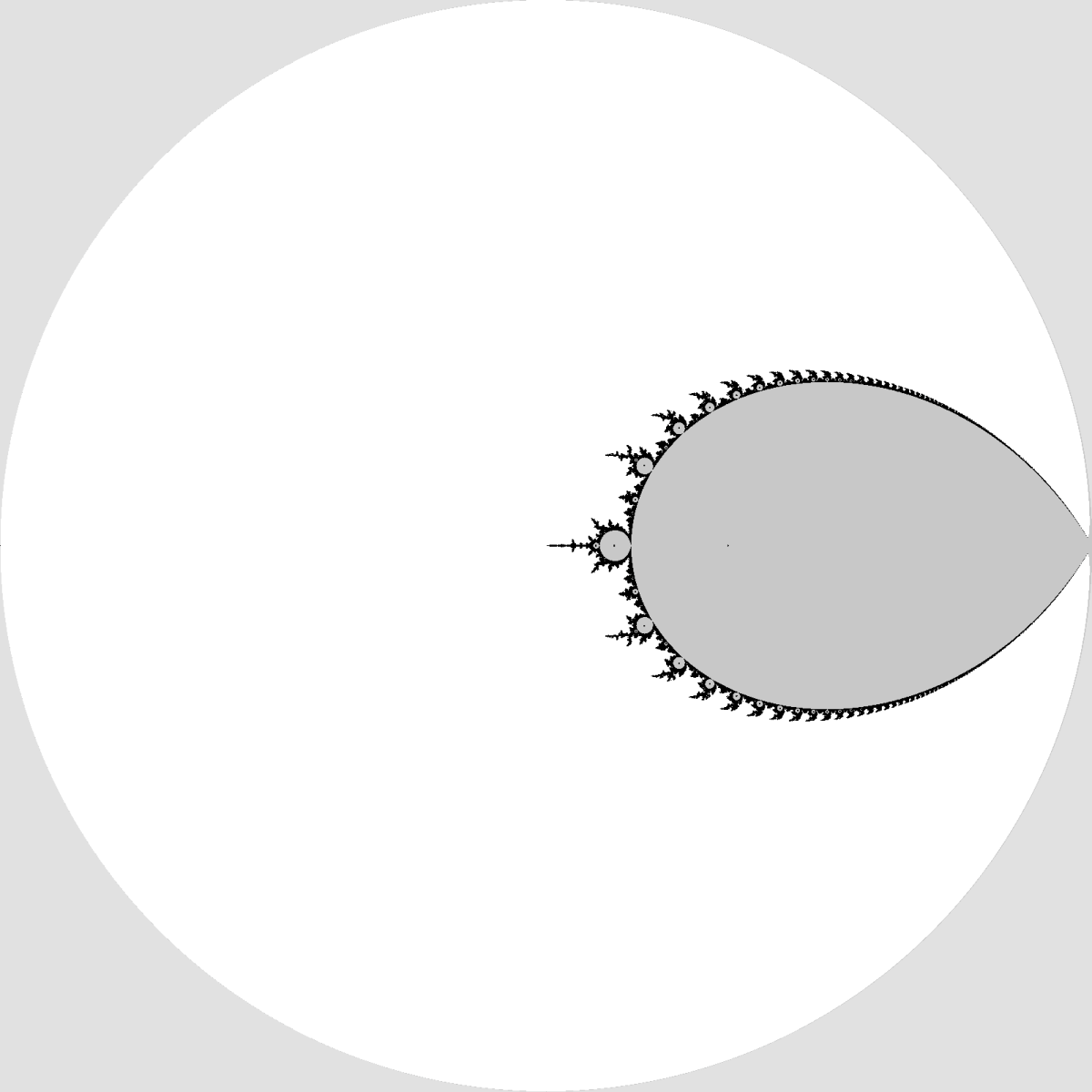}\hskip1.0cm\includegraphics{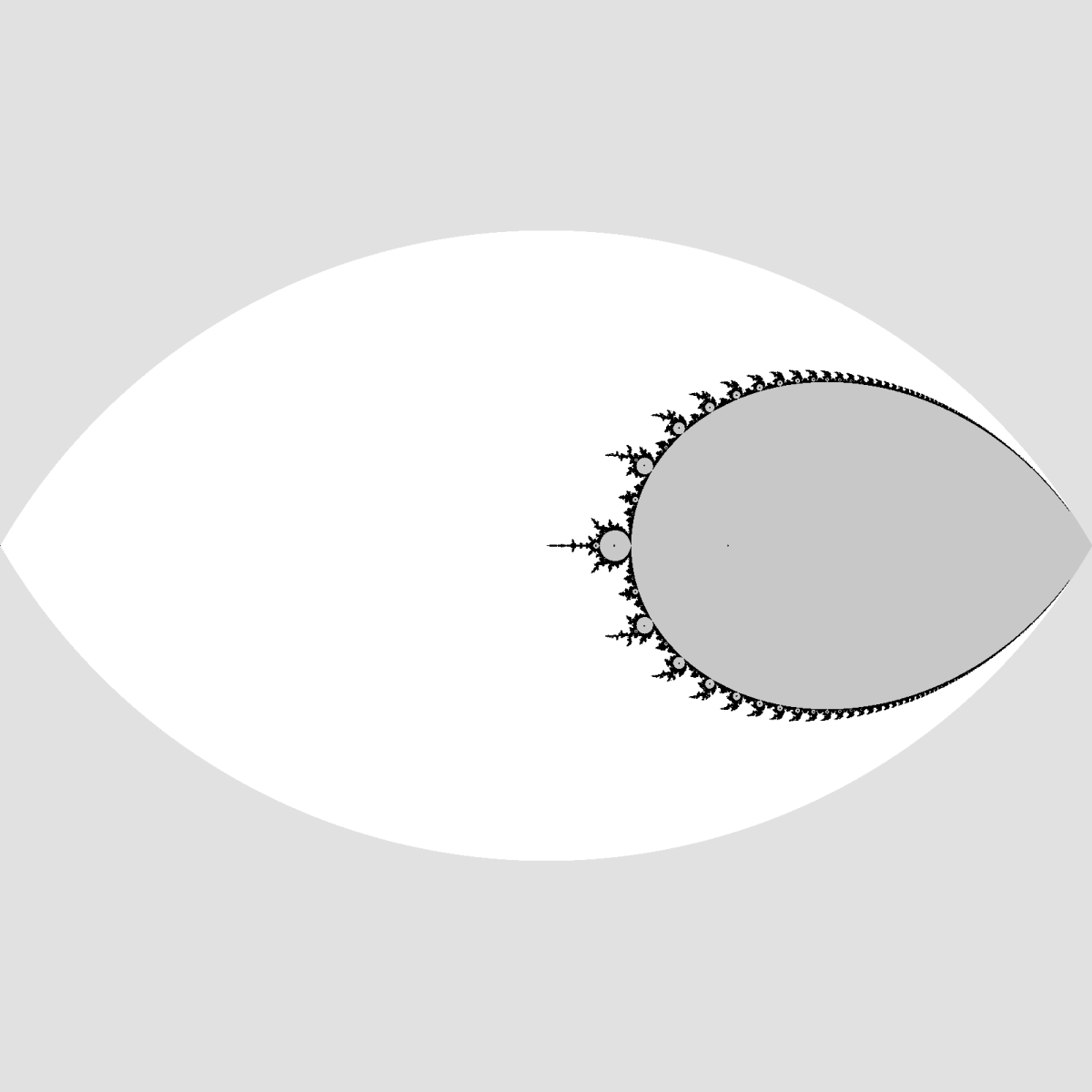}} 
\caption{\small The modular Mandelbrot set $\mathcal{M}_\Gamma$. On the left the white region is the round disc which has centre $a=4$ and radius $3$; on the right
it is the lune ${\mathcal L}_{\pi/3}$, of vertex angle $2\pi/3$. We prove that ${\mathcal M}_\Gamma\subset{\mathcal L}_{\theta}$ for some
$\theta$ with $\pi/3\le \theta<\pi/2$ (Theorem \ref{param-lune}).}
\label{mandelcorr1}
\end{center}
\end{figure}

An $(m:n)$ {\it holomorphic correspondence} on the Riemann sphere $\widehat \C$ is an $n$-valued function $F:z\to w$ (with $m$-valued inverse
$F^{-1}:w \to z$) defined implicitly by a polynomial equation $P(z,w)=0$, where $P$ has degree $m$ in $z$ and $n$ in $w$. 
In this paper we study the dynamics under iteration of members of a particular one complex parameter family of $(2:2)$ holomorphic correspondences, namely the family 
${\mathcal F}_a:z \to w$, given by  
$$\left(\frac{az+1}{z+1}\right)^2+\left(\frac{az+1}{z+1}\right)\left(\frac{aw-1}{w-1}\right)
+\left(\frac{aw-1}{w-1}\right)^2=3,$$ introduced by the first author, together with Christopher Penrose, in \cite{BP}.
For this family, there exists an open set $\mathcal K \subset \C$, called the  \textit{Klein Combination Locus}, such that for all $a \in \mathcal K$ 
the Riemann sphere admits a partition into two subsets completely invariant under $\F_a$, denoted $\Lambda_a$ and $\Omega_a$. The \textit{limit set} $\Lambda_a$ is closed and 
is the union of two closed subsets, the \textit{forwards limit set} $\Lambda_{a,-}$, restricted to which $\F_a$ is $(2:1)$, and the \textit{backwards limit set} $\Lambda_{a,+}$, 
on which $\F_a$ is $(1:2)$;
the intersection $\Lambda_{a,-}\cap\Lambda_{a,+}$ consists of a parabolic fixed point of $\F_a$ of multiplier $1$ (see Section \ref{prelim}). 
The complement $\Omega_a= {\widehat \C}\setminus \Lambda_a$ of the limit set is called 
the \textit{regular set}. The \textit{connectedness locus} $\mathcal C_\Gamma$ for the family $\F_a$ is the set of $a\in \C$ such that the limit set $\Lambda_a$ is connected (see Section \ref{prelim}). 
In \cite{BL1}, the current authors prove that when $a \in \mathcal C_\Gamma$ the correspondence $\F_a$ behaves like a parabolic quadratic map of the form 
$$P_A(z)=z+1/z+A$$ 
on a doubly pinched 
neighbourhood of $\Lambda_{a,-}$, using parabolic-like maps (see \cite{L}), and like the modular group on its complement. 
To make the statement precise, let us recall that, for every $A \in \C$, the map $P_A$ has a parabolic fixed point at $z=\infty$ with multiplier $1$, and basin of attraction $\A_A(\infty)$; 
the filled Julia set $K_A$ of $P_A$ is defined to be $K_A= \widehat \C \setminus \A_A(\infty)$. The Main Theorem in \cite{BL1} states that 
for every $a \in \mathcal C_\Gamma$, the correspondence $\F_a$ is hybrid conjugate on a doubly pinched neighbourhood of $\Lambda_{a,-}$ to a rational map of the form 
$P_A$ acting on a doubly pinched neighbourhood of its filled Julia set $K_A$, and on its regular set $\Omega_a$ the correspondence $\F_a$ is conformally conjugate 
to the pair of generators 
$$\alpha: z \to z+1\ \ \ \  {\rm and}\ \ \ \ \beta: z \to \frac{z}{z+1}$$ 
of the modular group $\Gamma=PSL(2,\Z)$ acting on the complex upper half-plane $\mathbb H$.
In other words, for  every $a \in {\mathcal C}_\Gamma$ the correspondence ${\mathcal F}_a$ is a \textit{mating} between a rational 
map of the form
$P_A:z \to z+1/z+A$ and $\Gamma$. 
In view of this theorem, we call the correspondences in the family $\F_a$ \textit{modular matings}.
We define the \textit{modular Mandelbrot set} $\m$ to be $\m:=\mathcal C_{\Gamma} \cap \overline \D(4,3)$. See Figure \ref{mandelcorr1}.\\

\begin{figure}
 
\begin{center}
\centering
\scalebox{0.8}{\includegraphics{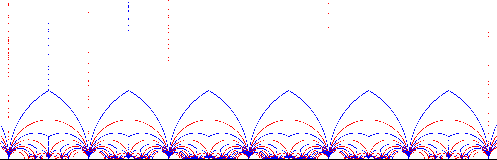}}

\vskip1.5cm
\centering
\scalebox{0.7}{\includegraphics[width= 7.0cm]{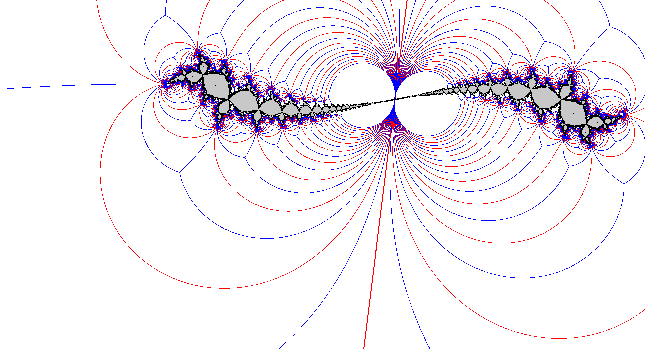}}
\end{center}
\caption{\small Above: tessellation of ${\mathbb H}$ by copies of a fundamental domain for
the modular group $\Gamma=PSL(2,{\mathbb Z})$. Below: mating between a `Douady rabbit' in 
$Per_1(1)$ and the modular group, realised by the correspondence ${\mathcal F}_a$ with $a=4.53926+0.439437i$. There are two copies of the rabbit, $\Lambda_{a,-}$ on the left
and $\Lambda_{a,+}$ on the right. They intersect in a parabolic fixed point $P$ for the correspondence (in the centre of the computer plot).
}\label{mating_pic}

\end{figure}

In the current paper we develop a complete dynamical theory for these modular matings.
For a general holomorphic correspondence on the Riemann sphere, analysing dynamics under iteration is a very challenging problem, even in the $(2:2)$ case. It is akin to the problem of analysing the dynamics of the group generated by a general pair of fractional linear transformations, and it is well known that for such a group it is a non-trivial problem
even to determine whether it is discrete. Isolated examples of correspondences (for example the arithmetic-geometric mean) are well understood. The family  of correspondences
$\F_a, a\in \mathcal{M}_{\Gamma}$, is the simplest $1$-parameter family to be fully explored, and as such it should provide a paradigm for the study of other families of correspondences displaying ``discreteness''  in their action on the sphere for some subset of parameter space.\\

\begin{figure}
\begin{center}
\scalebox{.5}{\includegraphics{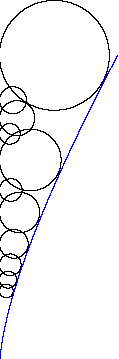}} \hspace{3cm} \scalebox{.5}{\includegraphics{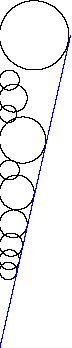}}
\caption{\small Discs in the $\log(\zeta)$-plane permitted by the Yoccoz inequality: on the left for matings between the modular group 
and quadratic polynomials, and on the right for the classical case of quadratic polynomials. In both diagrams the discs plotted correspond to all $0<p/q\le1/2$ with $q\le 8$.
}\label{Yoccoz_circles}
\end{center}
\end{figure}

For $a\in {\mathcal C}_\Gamma$, the Main Theorem in \cite{BL1} guarantees the existence of a canonical conformal homeomorphism 
$$\varphi_a:\widehat \C \setminus \Lambda({\mathcal F}_a) \to {\mathbb H},$$ conjugating the correspondence 
restricted to $\Omega({\mathcal F}_a)$ to 
the pair of generators $\alpha, \beta$ of the modular group defined above (see Figure \ref{mating_pic} and Section \ref{Boettcher}).
We call this map the \textit{B\"ottcher map} for the correspondence 
${\mathcal F}_a$, $a\in {\mathcal C}_\Gamma$, because it plays a role in our theory analogous to that of the classical B\"ottcher map 
$$\varphi_c:\widehat \C \setminus K(Q_c) \to \widehat \C\setminus \overline \D$$
which conjugates a quadratic polynomial $Q_c:z \to z^2+c$ on the complement of its filled Julia set $K(Q_c)$ to the map $z\to z^2$ on the complement of the closed unit disc, whenever $c$ lies in the Mandelbrot set ${\mathcal M}$, the set of $c\in \C$ such that $K_c$ is connected. 
Let us recall that external rays are the lines $R_\theta=\{z|\arg(\varphi_c(z))=\theta\}$,  preimages under the B\"ottcher map $\varphi_c$
of lines of constant argument in $\C$, and are important tools to understand the dynamics of a polynomial. Douady and Hubbard proved that {\it periodic rays land} (that is, there exists a limit, 
necessarily in $\partial K(Q_c)$, for 
$\varphi_c^{-1}(r e^{2\pi i \theta})$ as $r \to 1$), and that repelling periodic points are landing points of at least one and at most finitely many external rays. By considering properties of the linearisation of $Q_c$ in the neighbourhood of a repelling periodic cycle, Yoccoz proved, for $Q_c$ with $c\in {\mathcal M}$, an inequality constraining the (complex) logarithm of the multiplier of a repelling periodic cycle to lie in a certain disc in $\C$ (see \cite{H} and Figure \ref{Yoccoz_circles}), an inequality which he applied to estimate the size of the limbs of the Mandelbrot set.\\

There are strong parallels between the dynamical theories for the family $\F_a$ and the family $Q_c$, but also essential differences. The core similarity is that both theories are constructed using isometries of the hyperbolic metric on $\H$: the map $z\to z^2$ on 
the punctured disc $\C \setminus \overline \D$ lifts to the isometry $z\to 2z$ on its universal cover $\mathbb H$, and external rays in
$\C \setminus K(Q_c)$ are the preimages, under the lift of the B\"ottcher map, of geodesics in $\H$ from points 
of $\R\subset\partial\H$ to $\infty$. The core difference between the theories comes from fact that $z \to 2z$ is a {\it hyperbolic} isometry 
(it has two fixed points, $0$ and $\infty$ on $\partial \H$), whereas $z\to \alpha(z)$ and $z\to \beta(z$) are two {\it parabolic} isometries 
(with one fixed point each, $\infty$ and $0$, respectively).
The modular group endows
${\mathbb H}$ with a richer structure than that which the angle-doubling map gives to $\widehat \C\setminus \overline \D$,
but conversely there is no single point in ${\mathbb H}$ which has properties for $\Gamma$ analogous to those 
enjoyed by $\infty\in \widehat \C\setminus \overline \D$ for the angle-doubling map. For the family $\F_a$ the role of {\it external rays} in the Douady-Hubbard theory for $Q_c$ is replaced by that of {\it geodesics} in $\widehat{\mathbb C}\setminus\Lambda_a$ with initial point in $\Lambda_{a,-}$ and final point in $\Lambda_{a,+}$. The Douady-Hubbard theory
makes frequent use of a Green's function $G_c:{\mathbb C}\setminus K(Q_c) \to {\mathbb R}^{>0}$ which interacts with the dynamics
via $G_c(Q_c(z))=2G_c(z)$. The Green's function $G_a$ for $\Lambda_a$ has no such simple relationship with the two branches of $\F_a$,
but for points close to $\Lambda_{a,-}$ we prove that there exist lower and upper bounds for the multiplying factor $G_a(\F_a(z))/G_a(z)$
(Lemma \ref{alphabetapotential}, Section \ref{Green_fn}). This weaker property suffices for our purposes.
\\

In Proposition \ref{zero-accessible}  (Section \ref{Boettcher}) we prove  that
the inverse $\psi_a$ of the B\"ottcher map extends continuously to the end points of the imaginary axis
${\mathcal I}$ in ${\mathbb H}$ (the geodesic running from $0$ to $\infty$) sending both $0$ and $\infty$ to the parabolic fixed point 
$\Lambda_{a,-} \cap \Lambda_{a,+}$ (see Fig. \ref{mating_pic}), in other words that this geodesic {\it lands}, at both ends. 
In Proposition \ref{periodic_rays_land} (Section \ref{geod}) we prove that periodic geodesics land. 
Here, by a {\it periodic geodesic} we mean one which is
mapped to itself by some finite sequence of branches of $\F_a$: under the B\"ottcher map $\varphi_a$ it corresponds to a geodesic
in ${\mathbb H}$ fixed by some finite product of positive powers of $\alpha$ and $\beta$.
(Note that while ${\mathcal I}$ is not itself a periodic 
geodesic, its end-points are fixed by $\alpha$ and $\beta$ respectively.)
We then explain how one may code a geodesic in $\mathbb H$ from a point in $\mathbb R^-$ to a point in $\mathbb R^+$  by a bi-infinite sequence 
$S\in \{\alpha,\beta\}^\mathbb Z$
combining the continued fraction expansions of these two points. The $2$-sided shift operator on bi-infinite sequences imposes a dynamic on the space of geodesics from 
$\mathbb R^-$ to $\mathbb R^+$ in which the {\it periodic geodesics} play an organising role. 
By adapting a technique due to Benini and Lyubich \cite{BLyu}, 
in Section  \ref{repellor} we prove the deeper converse to Proposition \ref{periodic_rays_land} that 
\textit{repelling periodic points are landing points of 
periodic geodesics}. More precisely:

\begin{thm}\label{repelling_fp}

For every $a \in {\mathcal C}_\Gamma$, 

(i) every repelling fixed point of the $2$-to-$1$ restriction $f_a$ of ${\mathcal F}_a$ to a (doubly pinched) neighbourhood of 
$\Lambda_-({\mathcal F}_a)$ is the landing point of exactly one cycle of periodic geodesics;

(ii) the points of every repelling cycle of period $m>1$ of $f_a$ are landing points of at least one and at most two cycles of periodic geodesics.

\end{thm}

As a consequence of Theorem \ref{repelling_fp}(i), we may define
the {\it combinatorial rotation number} of a repelling fixed point $\hat{z}$ of $f_a$ (where $a \in {\mathcal C}_\Gamma$)
to be the rotation number of the unique periodic cycle of geodesics landing at $\hat{z}$. 
In Section \ref{Yoccoz_proof} we prove an inequality of Pommerenke-Levin-Yoccoz type for our correspondences
at such a fixed point:
\begin{thm} \label{Yoccoz}
If $a\in {\mathcal C}_\Gamma$ and $p_a$ is a repelling fixed point of $f_a$ which has
multiplier $\zeta$ and 
combinatorial rotation number $p/q\in {\mathbb Q}/{\mathbb Z}$, where $0<p/q<1$,
then $\log{\zeta}$
lies in the disc in the right hand half of ${\mathbb H}$ which has boundary 
circle tangent to the imaginary axis at $2\pi i p/q$ and radius $r_{p/q}$, where
$$r_{p/q}= \frac{2p\log(\lceil q/p\rceil+1)}{q^2} \quad\quad {\rm if} \quad 0< p/q \le 1/2,$$
and
$$r_{p/q}= \frac{2(q-p)\log(\lceil q/(q-p)\rceil+1)}{q^2} \quad\quad {\rm if} \quad 1/2\le p/q < 1.$$
\end{thm}

In the original Yoccoz inequality for quadratic polynomials, the corresponding radius 
$r_{p/q}$ is $1/q$ \cite{H}. Writing $2\pi i \nu$ for the imaginary part of $\log{\zeta}$, 
we deduce that whereas in the 
classical Yoccoz case for $c$ to be in $\mathcal M$ the value of $\log{\zeta}$ (where $\zeta$ is the derivative of $Q_c$ 
at a repelling fixed point) must lie in a strip whose width
goes to zero linearly with $\nu$, in our situation for $a$ to be in $\mathcal{M}_{\Gamma}$ the value of $\log{\zeta}$
has to lie in a strip whose width goes to zero at a rate proportional to $\nu^2\log(1/\nu)$
(for $\nu=1/q$ this is $(\log{q})/q^2$): see Figure \ref{Yoccoz_circles}. For repelling orbits of period greater than $1$ the bound 
on the modulus of the derivative is more complicated: see Remark \ref{higher_period} at the end of Section \ref{Yoccoz_proof} for an example.
\\

In Section \ref{param_lune} we apply Theorem \ref{Yoccoz} to prove that the modular Mandelbrot set $\mathcal{M}_{\Gamma}$
is contained in a subset of $\mathbb C$ of a particular form.
By the {\it lune} ${\mathcal L}_\theta$ in the parameter plane we shall mean the closed subset of the $a$-plane
bounded by the two arcs of circles which pass through the points $a=1$ and $a=7$ and meet the real axis
at angles $\pm\theta$ at these points (so the lune ${\mathcal L}_{\pi/2}$ is the disc $\overline \D (4,3)$).
The plot on the right in Figure \ref{mandelcorr1} suggests the possibility that ${\mathcal M}_\Gamma\subset {\mathcal L}_{\pi/3}$.
We do not know if this is the case, but we prove:

\begin{thm}\label{param-lune}
There exists an angle $\theta$ in the half-open interval  $\pi/3\le \theta <\pi/2$ such that 
$\mathcal{M}_{\Gamma}\subset {\mathcal L}_\theta$ and $\mathcal{M}_{\Gamma}$ only meets $\partial {\mathcal L}_\theta$ 
at the vertex $a=7$.
\end{thm}

Writing $L_{p/q}$ for the $p/q$-limb of $\mathcal{M}_{\Gamma}$, the set of parameter values where the fixed point of $\F_a$ is repelling or neutral of combinatorial rotation number
$p/q$, and writing $a_{p/q}$ for the root point of $L_{p/q}$, the value of $a$ where the derivative of $\F_a$ at its fixed point is $e^{2\pi i p/q}$, the computations
in the proof of Theorem \ref{param-lune} specialise to:

\begin{cor}\label{near_7}
As $p/q>0$ converges to zero, $a_{p/q}$ converges to $a=7$ tangentially to a straight line at angle $2\pi/3$ to the positive real axis and $diam(L_{p/q})$ converges to zero.
\end{cor}

Indeed the picture in the $a$-plane close to $a=7$ is the image of the lefthand plot in Figure \ref{Yoccoz_circles} under the transformation $\log{\zeta} \to 7+6(i\log{\zeta})^{2/3}$.\\

\begin{table}\label{compare}
\scriptsize{
\begin{tabular}{| c | c | c |}
\hline
{\it Quadratic polynomials $Q_c$} & {\it Quadratic correspondences ${\mathcal F}_a$} & {\it Comments and/or} \\
 & & {\it location of discussion}\\
\hline
\hline

$J(Q_c)$ quotient of ${\mathbb R}/{\mathbb Z}$ & $\partial\Lambda({\mathcal F}_a)$ quotient of 
$\widehat \R/(0\sim \infty)$ & Minkowski `?-function'\\

binary coding & continued fraction coding & relates these, \ref{Mink_qm} \\

\hline
$\varphi_c:\widehat \C\setminus K(Q_c) \cong \widehat \C\setminus \overline{\mathbb D}$
& $\varphi_a:\widehat \C \setminus\Lambda({\mathcal F}_a)\cong {\mathbb H}$ & `B\"ottcher conjugacy'\\

external rays  & external geodesics  & Sections \ref{Boettcher} and \ref{geod}\\

\hline
`periodic rays land' & `periodic geodesics land' & Section \ref{geod}, Prop. \ref{periodic_rays_land}\\

\hline

every repelling or parabolic &         every repelling  & Thm. \ref{repelling_fp} and
Section \ref{repellor},  \\

periodic point is the landing &       periodic point is the landing & (Hence $\mathcal{M}_{\Gamma}$, like $\mathcal M$, has\\

point of a periodic ray &  point of a periodic geodesic  & `no irrational limbs')\\

\hline
Yoccoz inequality  & Yoccoz inequality & Thm \ref{Yoccoz} and Section \ref{Yoccoz_proof}\\

$\sim 1/q$ & $\sim (\log{q})/q^2$ & \\

\hline
$M\subset \overline{D(0,2)}$ & $\mathcal{M}_{\Gamma}\subset {\mathcal L}_\theta$ & Thm \ref{param-lune} and Section \ref{param_lune}\\

\hline

$\Phi:\widehat \C \setminus M \cong \widehat \C\setminus\overline{\mathbb D}$ & $\Phi:U\setminus \mathcal{M}_{\Gamma}\cong V\subset{\mathbb H}$
& `Douady-Hubbard map'\\
parameter space rays & {tessellation of nbhd of $\mathcal{M}_{\Gamma}$} & Fig. \ref{mandelcorr} and \cite{BL4} \\
\hline
\end{tabular}
}

\caption{\it Similarities and differences (for $c \in \mathcal M$ and $a\in \mathcal{M}_{\Gamma}$)}
\end{table}

For the convenience of the reader, Table \ref{compare} exhibits parallels between the 
Douady-Hubbard theory for quadratic polynomials $Q_c$ and the theory developed below for the quadratic correspondences 
${\mathcal F}_a$.\\

In a sequel \cite{BL3} to the present paper we shall show that when $a$ lies in the parameter space lune ${\mathcal L}_\theta$
established by Theorem \ref{param-lune}, we may choose a dynamical space lune $L_a$, depending analytically on $a$, 
on which to perform a surgery
to transform the correspondence ${\mathcal F}_a$ into a rational map of the form $z \to z+1/z+A$, that is to 
say a quadratic rational map whose conjugacy class lies in Milnor's slice $Per_1(1)$. 
Then, using the theory of {\it holomorphic motions} we shall prove the following conjecture, in which
$\mathcal{M}_1$ denotes the connectedness locus of the family $Per_1(1)$:

\begin{conj} \label{homeo} There exists a homeomorphism $\chi:\mathcal{M}_{\Gamma}\to \mathcal{M}_1$  such that for each 
$a\in \mathcal{M}_{\Gamma}$ the correspondence ${\mathcal F}_a$ is a mating between the quadratic rational map 
$\chi(a)\in Per_1(1)$ and the modular group $\Gamma$.
\end{conj}

In \cite{PR} Carsten Petersen and Pascale Roesch prove that $\mathcal{M}_1$ is homeomorphic
to the classical Mandelbrot set $\mathcal M$. Thus proving Conjecture \ref{homeo} will 
yield a mathematical justification of the resemblance between $\mathcal{M}_{\Gamma}$ and $\mathcal M$ first noted experimentally 
almost three decades ago in \cite{BP}.\\

\noindent\textbf{The anti-holomorphic case.} 
Contemporaneously with our work, but independently of it, Lee, Lyubich, Makarov and Mukherjee \cite{LLMM} investigated matings between 
quadratic anti-rational maps and a discrete group generated by hyperbolic reflections, 
as part of their innovative programme of study of the dynamics  
of Schwarz reflection maps associated to quadrature domains. Their notion of a mating between a map and a group differs from ours, in that 
theirs is a single-valued function, which on an invariant closed connected subset has the behaviour of an anti-rational map on its filled Julia set, 
and which on the complement is associated to a discrete group of isometries of $\H$ generated by a reflection and a rotation. Despite the differences, there is a close 
relationship between the two settings. The dynamics of an (anti-holomorphic) Schwarz reflection maps in the one parameter family 
considered in \cite{LLMM} strongly parallels that of a $(2:1)$ restriction
of a holomorphic correspondence in the one parameter family investigated in \cite{BP}, \cite{BL1}, and the present paper.
See Remark \ref{compare} at the end of Section \ref{prelim},
where the relationship between the two theories is elucidated. \\

{\bf Acknowledgement.} 
This work was supported by the Serrapilheira Institute (grant number Serra-1811-26166), the
EU Marie-Curie IRSES Brazilian-European partnership in Dynamical Systems (FP7-PEOPLE-2012-IRSES 318999 BREUDS), the Funda\c{c}\~ao de amparo a pesquisa do estado de S\~ao Paulo (Fapesp, process 2013/20480-7),
the Cnpq (406575/2016-19), and the prize L'ORÉAL-UNESCO-ABC Para Mulheres na Ciência.

\section{Preliminaries: limit sets for $\F_a$}\label{prelim}
In this Section we note some general definitions and properties of $(2:2)$ holomorphic correspondences,
and recall terminology and results from \cite{BL1}, which we shall use throughout the 
paper. \\  

A $(2:2)$ {\it holomorphic correspondence} on the Riemann sphere is a $2$-valued function $F:z\to w$ (with $2$-valued inverse
$F^{-1}:w \to z$) defined implicitly by a polynomial equation $P(z,w)=0$, where $P$ has the form
$$P(z,w)=(az^2+bz+c)w^2+(dz^2+ez+f)w+(gz^2+hz+j).$$
The graph of $F$, the Riemann surface $P(z,w)=0$, supports two natural involutions: $I_-$ which interchanges the two values of $z$ which map to the same $w$, 
and $I_+$ which interchanges the two values of $w$ which are images of the same $z$. 
We say
that $F$ is a {\it map of triples} if $(I_+I_-)^3$ is the identity. Equivalently the immediate images $w_2,w_3$ of any point $z_1$,
and the immediate pre-images of $w_2$ and $w_3$ fit together in a closed diagram:

\noindent
\begin{picture}(270,100)

\put(150,20){$z_3$}
\put(150,44){$z_2$}
\put(150,68){$z_1$}
\put(184,20){$w_3$}
\put(184,44){$w_2$}
\put(184,68){$w_1$}
\put(162,22){\vector(1,1){20}}
\put(162,26){\vector(1,2){20}}
\put(162,42){\vector(1,-1){20}}
\put(162,50){\vector(1,1){20}}
\put(162,66){\vector(1,-2){20}}
\put(162,70){\vector(1,-1){20}}

\end{picture}

which at isolated points may take one of the degenerate forms below:

\noindent
\begin{picture}(320,70)

\put(100,20){$z_2$}
\put(100,47){$z_1$}
\put(135,20){$w_2$}
\put(135,47){$w_1$}
\put(112,22){\vector(1,0){20}}
\put(112,25){\vector(1,1){20}}
\put(112,46){\vector(1,-1){20}}
\put(200,34){$z_1$}
\put(235,34){$w_1$}
\put(212,36){\vector(1,0){20}}

\end{picture}

In the diagram on the left, $(z_2,w_1)$ is a fixed point of $I_-$ and $(z_1,w_2)$ is a fixed point of $I_+$. 
We remark that $dw/dz=0$ for the forwards branch $z_2 \to w_1$, and
$dz/dw=0$ for the backwards branch $w_2 \to z_1$; thus the four points in the diagram are all  
forwards or backwards critical points or values of the correspondence. In the diagram on the right, $(z_1,w_1)$ 
is fixed by both $I_-$ and $I_+$; the point $(z_1,w_1)$ is a double point of the Riemann surface $P(z,w)=0$.\\

It follows from the definition, or equivalently from the diagrams above, that every map of triples necessarily factorises into the 
{\it deleted covering correspondence} $Cov_0^C$ of  a cubic 
$C:\widehat \C \to \widehat \C$ post-composed by a M\"obius transformation $M$.
Here $C$ is the map which identifies each triple of points
$\{z_1,z_2,z_3\}$ of the diagram condition to a single point, and $M$ is the map which sends each $z_j$ in the diagram to 
the corresponding $w_j$.
Up to pre- and post-multiplication by M\"obius transformations, $C$ can take one of three forms, depending on the numbers 
and types of critical points of the correspondence:

\begin{enumerate}[label=(\roman*)]
\item
$C(z)=z^3$;
\item
$C(z)=z^3-3z$;
\item
$C(z)=z^2(z+a)/(z+1)\ (a\ne 0,1,9)$.
\end{enumerate}

\begin{figure}
 
\begin{center}
\begin{tabular}{rcl}
\begin{minipage}{4cm}
\centering
 \psfrag{a}{\tiny $a$}
 \psfrag{A}{\tiny $ \Delta_J^{st}$}
  \psfrag{-2}{\tiny $-2$}
 \psfrag{1}{\tiny $1$}
  \psfrag{2}{\tiny $2$}
 \psfrag{C}{\tiny $\Delta_{Cov}^{st}$}
 \psfrag{L}{}
\includegraphics[width= 4cm]{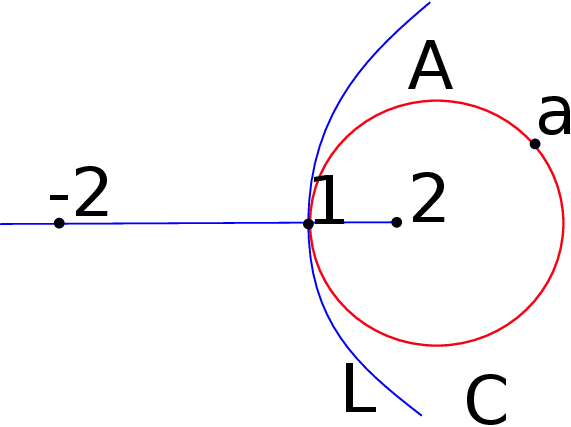}

 \end{minipage}
 
 \hspace{2.0cm}
\begin{minipage}{4cm}
\centering
 \psfrag{a}{\tiny $a$}
 \psfrag{A}{\tiny $ \Delta_J^{st}$}
  \psfrag{F}{\tiny $\F_a^{-1}(\Delta_J^{st})$}
 \psfrag{e}{\tiny $\F_a^{-2}(\Delta_J^{st})$}
  \psfrag{2}{\tiny $(2:1)$}
 \psfrag{1}{\tiny $(1:2)$}
  \psfrag{3}{\tiny $(1:1)$}
 \centering
\includegraphics[width= 4cm]{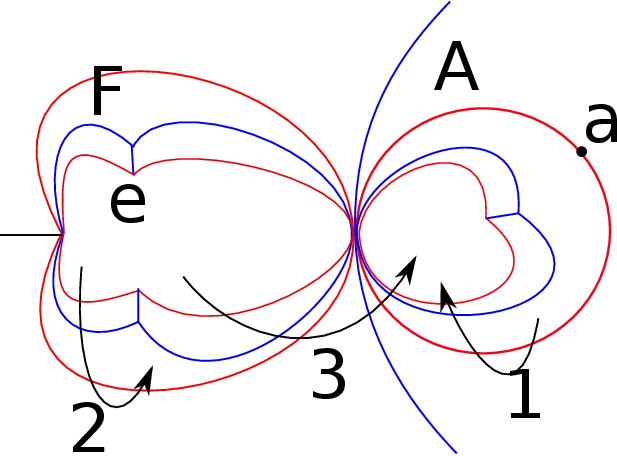}\end{minipage}

\end{tabular}
\end{center}
\caption{\small On the left: fundamental domains for $\F_a$ in $Z$ coordinates. On the right: first iterations of $\F_a$
}\label{domZ}

\end{figure}

These forms correspond to the three possible topological types for the graph of a $(2:2)$ holomorphic correspondence, namely two (intersecting) spheres, a single (self-intersecting) sphere, or (the generic case) an elliptic curve.
Our interest is in the second case, and the particular sub-case that the M\"obius transformation $M$ is an involution (which we 
denote $J$). This restriction reduces us to the two (complex) parameter family of {\it reversible maps of triples} which have graph
a single sphere. Our final restriction, reducing us to  one parameter, is to require that one of the fixed points of $J$ be also a critical point of $Cov_0^C$. 
Hence, we can write $\F_a=J_a\circ Cov_0^C$, where $Cov^C_0$ denotes the deleted covering correspondence of $C(z)=z^3-3z$, that is to say the $2$-to-$2$ correspondence $Z \to W$ on $\widehat \C$
defined by $$\frac{C(Z)-C(W)}{Z-W}=0, \quad {\rm that\ is,} \quad Z^2+ZW+W^2=3,$$
and $J_a$ denotes the conformal involution on $\widehat \C$ which has fixed points $Z=1$ and $Z=a$, namely $$J_a(Z)=\frac{(a+1)Z-2a}{2Z-(a+1)}.$$
Thus  the correspondence $\F_a: Z \rightarrow W$ is given by
$$Z^2+ZJ(W)+(J(W))^2=3,$$
which under the change of coordinates  
$$Z=(az+1)/(z+1), \mbox{ and } W=(aw+1)/(w+1)$$
is the relation
$$\left(\frac{az+1}{z+1}\right)^2+\left(\frac{az+1}{z+1}\right)\left(\frac{aw-1}{w-1}\right)
+\left(\frac{aw-1}{w-1}\right)^2=3.$$
By construction, the involution $J_a$ conjugates ${\mathcal F}_a$ to its inverse ${\mathcal F}_a^{-1}$. In the $z$-coordinate $J$ is the involution $z \to -z$ 
(whatever the value of  $a$). The {\it forwards critical points} of ${\mathcal F}_a$ (the points where $dW/dZ$ vanishes for a branch $Z\to W$ of ${\mathcal F}_a$, see \cite{BP} Section 2) are $Z=-1$ and $+1$. The corresponding {\it forwards critical values} are $W=J(+2)$ and $J(-2)$ respectively. The point $Z=\infty$ is a {\it double point} of ${\mathcal F_a}$ (see \cite{BP} Section 2).\\

As $\F_a=J_a \circ Cov_0^C$, the two images of $Z \in \C$ under $\F_a$ are the images under $J_a$ of the two preimages of 
$C(Z)=Z^3-3Z$ in $\C$ different from $Z$.
Using this decomposition, we will now construct fundamental domains for the action of $\F_a$. 
As $Z=1$ is a critical point of the cubic $C$, and $C(-2)=C(1)$, we deduce that $Cov_0^C(-2)=1$ under both branches, and $Cov_0^C((-\infty,-2])$ is  
a smooth curve running through $Z=1$, symmetric about the real axis, and asymptotic to the directions $e^{\pm\pi i/3}$ as $Z$  tends to $\infty$.
On the other hand, $J_a$ fixes the circle through  $Z=1$ and $Z=a$ having centre
on the real axis, and sends points inside the disc bounded by this circle to the exterior of the disc, and vice versa.
The {\it standard fundamental domain} $\Delta^{st}_{Cov}$ is the subset of the $Z$-plane to the right of the curve $Cov_0^C((-\infty,-2])$.
The  {\it standard fundamental domain} $\Delta^{st}_J$ is the subset of the $Z$-plane exterior to the circle through $Z=1$ and $Z=a$ which has centre
on the real axis (see Figure  \ref{domZ} on the left). Define $\mathcal D:=\{a:|a-4|\le3\}$, and note that for $a=1$ the correspondence is not defined. When $a \in \mathcal{D}\setminus \{1\}$, we have
$$\Delta^{st}_{Cov} \cup  \Delta^{st}_J =\widehat \C \setminus \{1\}.$$
Moreover, when $|a-4|<3$ these satisfy the condition that their
boundaries are transverse to the parabolic axis at $Z=1$ (see Proposition 3.5 in \cite{BL1} for a proof). When $a=7$ there are three attracting-repelling directions at $Z=1$ and the domain boundaries
are transverse to all three.\\

More generally, the {\it Klein Combination Locus}, $\mathcal K$, is the set of values of the parameter $a \in \widehat \C$ such that 
there exist fundamental domains $\Delta_{Cov}$ for $Cov_0^C$ and $\Delta_J$ for $J$ which together
cover all of $\widehat \C$ except for the single point $Z=1$ ($z=0$) (in the present article, as in \cite{BL1},
fundamental domains are open, so do not include their boundaries). A pair of fundamental domains
$(\Delta_{Cov}, \Delta_J)$ satisfying this condition is called a {\it Klein combination pair}. 
Since the standard fundamental domains are a Klein combination pair, we have that $\mathcal D \setminus \{1\} \subset \mathcal K$.
When $a \in {\mathcal K}$ we can always choose $\Delta_{Cov}$ and $\Delta_J$ such that the Jordan curves 
$\partial\Delta_{Cov}$ and $\partial\Delta_J$ are smooth at the parabolic fixed point $Z=1$ ($z=0$) of ${\mathcal F}_a$, and 
transverse to the attracting-repelling axis there (see Proposition 3.8 in \cite{BL1}).\\

\begin{figure}
 
\begin{center}
\begin{tabular}{rcl}
\begin{minipage}{3cm}
\centering
\includegraphics[width= 3.0cm]{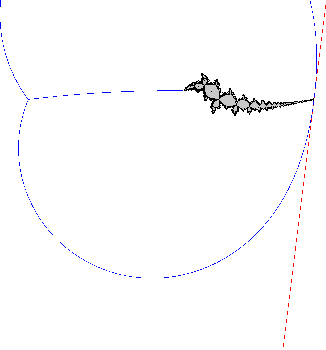}

 \end{minipage}
 
 \hspace{1.0cm}
\begin{minipage}{3cm}
\centering
\includegraphics[width= 3.0cm]{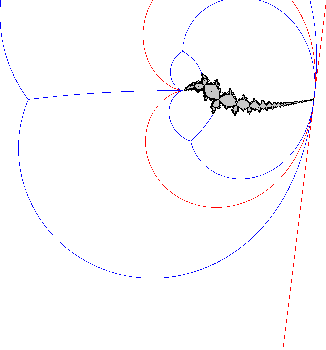}
\end{minipage}

 \hspace{1.0cm}
\begin{minipage}{3cm}
\centering
\includegraphics[width= 3.0cm]{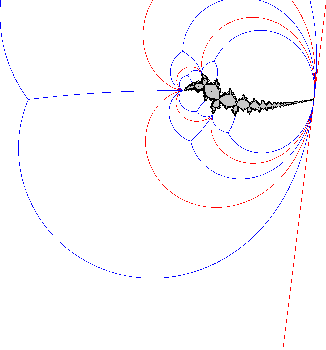}
\end{minipage}

\end{tabular}
\end{center}
\caption{\small Boundaries of the standard domains, plotted in the $z$-coordinate for the correspondence ${\mathcal F}_a$ with $a=4.53926+0.439437i$, together with  their images 
under ${\mathcal F}_a^{-1}$ and ${\mathcal F}_a^{-2}$. The domain $\Delta_{Cov}$ is the 
part of $\widehat \C$ outside the outer closed blue curve, $\Delta_J$ is the part of 
$\widehat \C$ to the left of the slanting straight red line. 
The set shaded grey is the backwards limit set 
$\Lambda_{a,-}=\bigcap_{n=0}^\infty {\mathcal F}_a^{-n}(\bar{\Delta}_J)
=\bigcap_{n=0}^\infty {\mathcal F}_a^{-n}(\widehat \C\setminus \bar{\Delta}_{Cov})$.
}\label{stand_domains}

\end{figure}

By construction, when $a\in {\mathcal K}$ and $(\Delta_{Cov}, \Delta_J)$ is a Klein combination pair, ${\mathcal F}_a$ restricted to the domain $\Delta_{Cov}$
is a $1$-to-$2$ correspondence  (see Figure  \ref{domZ} on the right, and Proposition 3.4 in \cite{BL1} for a proof): 
$${\mathcal F}_a:\Delta_{Cov} \to J(\widehat \C\setminus \Delta_{Cov})\subset \Delta_{Cov}.$$
Provided that $\Delta_{Cov}$ and $\Delta_J$ have smooth boundaries at the parabolic point and that these are transverse to the attracting-repelling axis there,
$$\Lambda_{a,+}:=\bigcap_{n=0}^\infty {\mathcal F}_a^n(\Delta_{Cov})$$
is independent of the choice of the Klein combination pair $(\Delta_{Cov}, \Delta_J)$. We call $\Lambda_{a,+}$ the {\it forwards limit set} of ${\mathcal F}_a$.
Moreover, $\widehat \C\setminus\Delta_{Cov} \subset J(\Delta_{Cov})$ and the restriction of ${\mathcal F}_a$ is a 
$2$-to-$1$ map from the first of these sets onto the second. We define the {\it backwards limit set} of 
${\mathcal F}_a$ (see Figure \ref{stand_domains}) to be
$$\Lambda_{a,-}:=
\bigcap_{n=0}^\infty {\mathcal F}_a^{-n}(\widehat \C\setminus\Delta_{Cov}).$$
Also by construction, $\Lambda_{a,-}=J(\Lambda_{a,+})$ and $\Lambda_{a,-}\cap\Lambda_{a,+}$ consists of a single point, the parabolic fixed point ($Z=1$) of 
${\mathcal F}_a$. We denote
$\Lambda_{a,-} \cup \Lambda_{a,+}$ by $\Lambda_a$.
Finally, ${\mathcal F}_a$ acts properly discontinuously on $\Omega_a=\widehat \C\setminus\Lambda_a$, with fundamental domain 
$\Delta_{corr}:=\Delta_{Cov} \cap \Delta_J$.\\

Further restricting the codomain of the $2$-to-$1$ branch of ${\mathcal F}_a$ defined above to 
$\overline{\Delta}_J$,  yields a $2$-to-$1$ map which we denote $f_a$. This inherits the property that $f_a^{-1}(\overline{\Delta}_J)\subset\overline{\Delta}_J$. Modulo boundary curves,
$\overline{\Delta}_J\setminus f_a^{-1}(\overline{\Delta}_J)$ is 
the disjoint union of $\Delta_{Cov}\cap\Delta_J$ and its two images under $Cov_0^C$. 
The map $f_a$ has unique critical point $Z=-1$ in $\Delta_J$ (the other forwards critical point $Z=+1$ of 
${\mathcal F}_a$ being a fixed point on the boundary $\partial\Delta_J$).\\

The connectedness locus ${\mathcal C}_\Gamma$ is the set of parameter values $a\in {\mathcal K}$ such that $\Lambda_{a,-}$ (or equivalently
$\Lambda_{a,+}$ or $\Lambda_a$) is connected. 
The {\it modular Mandelbrot set}, $\mathcal{M}_{\Gamma}$, is ${\mathcal C}_\Gamma \cap \mathcal{D}$.
It is proved in \cite{BL3} that 
$\mathcal{M}_{\Gamma}={\mathcal C}_\Gamma$.

 \begin{remark}
Replacing the involution $J_a$ in the construction above by reflection $\iota_a$ in the circle which has centre $Z=a$ and passes through $Z=1$, yields a family of $(2:2)$ 
anti-holomorphic correspondences 
$${\mathcal G}_a:=\iota_a\circ Cov_0^C,$$
where, as in the definition of $\F_a$, the function $C$ is the cubic $C(Z)=Z^3-3Z$. Let $\Delta_a$ denote the disc in the $Z$-plane bounded by the circle, and (by analogy with 
the definition of $\Delta_J$) let $\Delta_\iota$ denote its complement in $\hat{\C}$. The Klein combination condition becomes the condition that there exists a transversal 
$\Delta_{Cov}$ for $C$ which contains $\overline{\Delta}_a \setminus \{Z=1\}$, or, equivalently, that $C$ is injective on $\overline{\Delta}_a$. 
The correspondences in the family ${\mathcal G}_a$ which satisfy 
 this condition are intimately related to the family of (anti-holomorphic) Schwarz reflection maps $\sigma_a$ investigated in \cite{LLMM}.
This relationship may be deduced from the discussion in Section 10 of \cite{LLMM}, but to assist readers we shall describe it explicitly here.\\
Before we get into the details, recall that in \cite{LLMM} the notion of a {\it mating} between an anti-rational map and a group is not a correspondence, but a {\it map}, which on an 
invariant simply-connected closed subset of the Riemann sphere is conjugate to an anti-rational map on its filled Julia set, and on the complement of this subset behaves 
like a certain {\it map} associated to a group of automorphisms of $\H$.
Now let $\Sigma$ denote the quotient sphere $\hat{\C}/C$, and suppose that $a$ is such that $C$ is injective on the closure $\overline{\Delta}_a$ of $\Delta_a$. 
Let $g_a:{\mathcal G}_a^{-1}(\Delta_\iota) \to \Delta_\iota$ denote the pinched anti-quadratic-like $2$-to-$1$ branch of ${\mathcal G}_a$ restricted to $\Delta_\iota$ as domain and codomain, 
in analogy to the branch $f_a$ of $\F_a$ defined a few lines before this Remark, and note that $\iota_a\circ g_a\circ \iota_a$ is the corresponding $2$-to-$1$ restriction of ${\mathcal G}_a^{-1}$
to $\Delta_a=\iota_a(\Delta_\iota)$. 
As $C$ is univalent on $\Delta_a$, its image $C(\Delta_a)\subset \Sigma$ is the {\it quadrature domain} $\Omega_a$ defined in \cite{LLMM}, Section 3.2. When restricted to
$\Delta_a$, the function $C$ conjugates $\iota_a\circ g_a\circ \iota_a$ to a $2$-to-$1$ map from a subset of $\Omega_a$ onto the whole of $\Omega_a$. 
This $2$-to-$1$ map is {\it precisely} the Schwarz reflection $\sigma_a$ (\cite{LLMM}, Section 3.2) associated to the quadrature domain $\Omega_a$, as 
is apparent from the formula in Proposition 2.3 of \cite{LLMM} for the Schwarz reflection map associated to a quadrature domain in general.
We deduce that the restriction of $C\circ \iota_a$ to $\Delta_\iota$ conjugates the pinched anti-quadratic-like map $g_a$ on $\Delta_\iota$ to the pinched anti-quadratic-like 
map $\sigma_a$ on $\Omega_a$. Under this conjugacy the backwards limit set $\Lambda_{a,-}=\bigcap_{n=0}^\infty g_a^{-n}(\Delta_\iota)$ of ${\mathcal G}_a$ is carried 
to the non-escaping set $K_a=\bigcap_{n=0}^\infty \sigma_a^{-n}(\Omega_a)$ of $\sigma_a$, and the tiling on $\Delta_\iota \setminus \Lambda_{a,-}$ is carried to that on 
$\Omega_a\setminus K_a$. Thus for any value of $a$ such that the correspondence ${\mathcal G}_a$ is a mating in the sense of \cite{BP}, \cite{BL1} and the present paper, 
between a quadratic anti-rational map and the group of automorphisms of $\H$ obtained from the modular group by substituting $z\to 1/\bar{z}$ for $z\to -1/z$, the Schwarz reflection 
map $\sigma_a$ is a mating  between the same map and group in the sense of \cite{LLMM}, and vice versa.
\end{remark}

\section{The B\"ottcher map}\label{Boettcher}

By Theorem A of \cite{BL1}, for every $a \in \mathcal{C}_\Gamma$ there exists a conformal homeomorphism 
$\varphi_a:\Omega_a \to {\mathbb H}$ which conjugates the two branches of ${\mathcal F}_a |_{\Omega_a}$ to  
the automorphisms $\alpha:z \to z+1$ and $\beta:z \to z/(z+1)$ of ${\mathbb  H}$. This $\varphi_a$ is unique, since any
automorphism $h$ of ${\mathbb H}$ which conjugates both $\alpha$ to itself and $\beta$ to itself is necessarily the identity
(to see this, observe that $h$ must fix both $\infty$ and $0$, and that therefore $h$ has the form $h(z)=\lambda z$; 
the fact that $h^{-1}\alpha h(z)=\alpha(z)$ for all $z\in {\mathbb H}$ now implies that $\lambda=1$). 
We shall refer to $\varphi_a$ as the `B\"ottcher map'. It is analogous to the map:
$$\varphi_c:{\mathbb C}\setminus K(Q_c)\to {\mathbb C}\setminus\overline{\mathbb D}$$ defined by Douady and Hubbard
\cite{DH}, from the complement of the filled Julia set $K(Q_c)$ 
of a quadratic polynomial $Q_c(z)=z^2+c$ (with $c\in \mathcal M$) to the complement of the closed unit disc, conjugating $Q_c$ to $z\to z^2$. \\

{\bf Notation.}
For $a\in {\mathcal C}_\Gamma$, the two branches of the correspondence ${\mathcal F}_a$ become 
(single-valued) homeomorphisms when we restrict ${\mathcal F}_a$ to $\Omega=\Omega({\mathcal F}_a)$. 
We denote these homeomorphisms by $g:\Omega \to \Omega$ and $h:\Omega \to \Omega$, where $g$ 
corresponds under $\varphi_a$ to $\alpha$, and $h$ corresponds to $\beta$. \\

We recall that $PSL(2,{\mathbb Z})$  is the free product of the subgroups $C_2$ generated by $\sigma:z \to -1/z$
and $C_3$ generated by $\rho: z \to -1-1/z$.  Under the B\"ottcher map $\varphi_a$ the branch $g$ of 
${\mathcal F}_a|_{\Omega_a}$ corresponds to 
 $\sigma\rho$ and $h$ corresponds to  $\sigma\rho^{-1}$. Thus $g^{-1}hg^{-1}$ corresponds to $\sigma$ and 
$g^{-1}h$ corresponds to $\rho$. We deduce:

\begin{lemma}\label{Boettcher-prop}
For $a \in {\mathcal C}_\Gamma$, the B\"ottcher map $\varphi_a$ is the Riemann mapping of the simply-connected open set $\Omega({\mathcal F}_a)$,
$$\varphi_a: \Omega({\mathcal F}_a)\to {\mathbb H},$$  
normalised to send the fixed point of 
$g^{-1} h g^{-1}$ to $i \in {\mathbb H}$ and the fixed point of $g^{-1}h$ to 
$(-1+i\sqrt{3})/2\in {\mathbb H}$.
\end{lemma}

Pursuing the same route as that followed by Douady and Hubbard for polynomials, we now consider the question as to 
how far the inverse of $\varphi_a$ extends to a continuous map from
the boundary of ${\mathbb H}$ to the boundary of $\Omega=\Omega({\mathcal F}_a)$.

\begin{prop}\label{zero-accessible}
The inverse 
$\psi_a:{\mathbb H}\to \Omega$
of the Riemann mapping $\varphi_a$ extends continuously to 
$0$ and $\infty \in \widehat \R=\partial {\mathbb H}$, sending both these points to the fixed point $Z=1$ 
of ${\mathcal F}_a$. 
\end{prop}

{\bf Proof} Let $\mathcal I$ denote the imaginary axis in ${\mathbb H}$, the geodesic which runs from 
$0\in \partial {\mathbb H}$ to $\infty\in \partial {\mathbb H}$.  The homeomorphism 
$\psi_a:{\mathbb H}\to \Omega$ sends $\mathcal I$ and its
images under the cyclic group generated by $\alpha:{\mathbb H} \to {\mathbb H}$ to an arc $\psi(\mathcal I)$ in $\Omega$
and its images under the cyclic group generated by $g:\Omega \to \Omega$. 
The points $it \in \mathcal I$, $t \to \infty$, lie on horocycles of $\alpha$ in a family converging
to the fixed point $z=\infty$ of $\alpha$, on the boundary of ${\mathbb H}$. So their images $\psi_a(it)$, $t \to \infty$,
lie on horocycles of $g$ in a family which converges to the fixed point $Z=1$ of $g$ (see Step 3 
in the proof of Theorem A in \cite {BL1}). Thus 
the boundary point $Z=1$ is accessible from within $\Omega$ by the path $\psi_a(\mathcal I)$, and 
setting $\psi_a(\infty)=1$ therefore extends $\psi_a$ continuously to the end point $z=\infty$ of $\mathcal I$. 
The proof for the other end of $\mathcal I$ is similar, with $\beta$ in place of $\alpha$ and $h$ in place of $g$.
\qed

\begin{cor} The inverse $\psi_a$ of $\varphi_a$ extends continuously to every $p/q \in{\mathbb Q}\subset\widehat \R$.
\end{cor}

{\bf Proof} This follows immediately from Proposition \ref{zero-accessible}, since the orbit of $0\in {\mathbb R}$
under $PSL(2,{\mathbb Z})$ is ${\mathbb Q}$. \qed
\\

We remark that the rationals perform the same role here as that of the {\it dyadic} rationals (those with finite binary expressions) 
in the case of quadratic polynomials. In Subsection \ref{per-geo} the {\it quadratic irrationals} will come into the picture for
$PSL(2,{\mathbb Z})$, playing a role analogous to that played by the non-dyadic rationals for quadratic polynomials.

\subsection{Digression: the Douady-Hubbard map and a tessellation of a neighbourhood of $\mathcal{M}_{\Gamma}$}\label{digression}

Recall that for quadratic polynomials $Q_c:z \to z^2+c$, Douady and Hubbard constructed a 
canonical bijection 
$\Phi:{\mathbb C}\setminus \mathcal M\to {\mathbb C}\setminus \overline{\mathbb D}$ by means of the ingenious 
assigment $\Phi(c)=\varphi_c(c)$. This map $\Phi$ enabled them to investigate $\mathcal M$ via {\it parameter rays} 
in ${\mathbb C}\setminus \mathcal M$. It is natural to ask whether there is an analogous construction in the parameter space of our family of correspondences. There is, and we outline it here, but we postpone the
proof to a separate article \cite{BL4}, as this will be easier with the aid of some of the methods and results of \cite{BL3}.
\\

\begin{figure}
\begin{center}
\scalebox{0.38}{\includegraphics{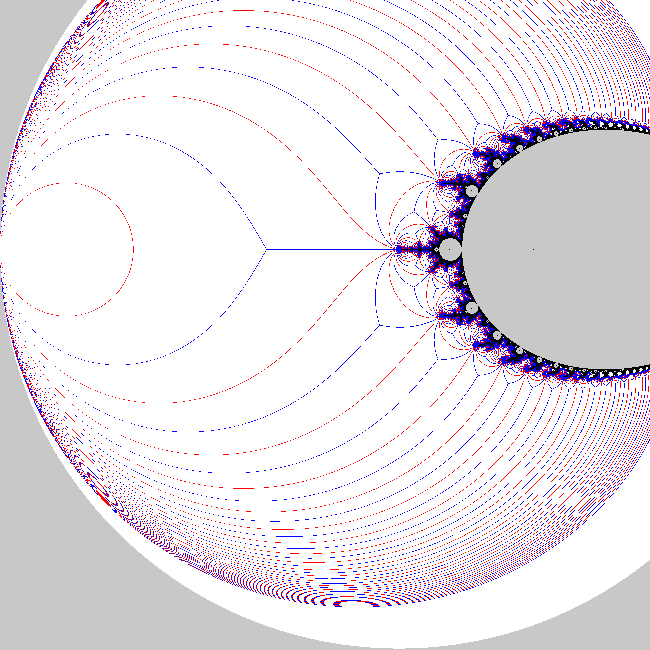}} 
\hskip1.5cm
\scalebox{.38}{\includegraphics{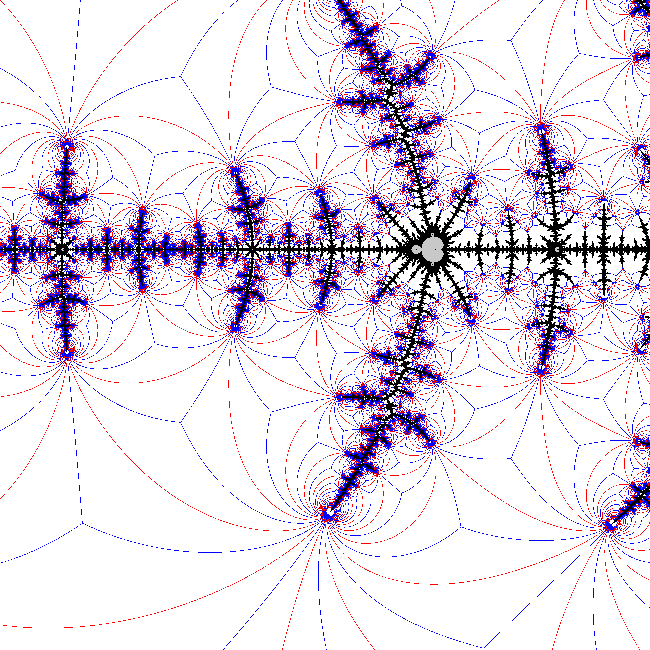}} 
\caption{\small Plot of the modular Mandelbrot set $\mathcal{M}_{\Gamma}$. The round white disc ${\mathcal D}$ has centre $a=4$ and radius $3$:
the tessellation on ${\mathcal D}\setminus \mathcal{M}_{\Gamma}$ is 
discussed in Section \ref{digression}. The zoom on the right has at its centre the small 
copy of $\mathcal{M}_{\Gamma}$ on the real axis which corresponds to period $3$.}\label{mandelcorr}
\end{center}
\end{figure}

B\"ottcher's conjugacy $\varphi_c$
exists between every quadratic polynomial $Q_c$ and the squaring map $z\to z^2$ on a neighbourhood of $\infty$, whether or not $c \in \mathcal M$. Douady and Hubbard observed that this conjugacy can always be extended to an open set containing the critical value $c$ of $Q_c$, although it can only be extended to the whole
of ${\mathbb C}\setminus K(Q_c)$ if $c \in \mathcal M$.  For a correspondence ${\mathcal F}_a$, we have no
canonical point analogous to $\infty$ from which to begin: instead we start from a Klein combination pair $(\Delta_{Cov},\Delta_J)$ for ${\mathcal F}_a$. For every $a$ in the Klein combination locus ${\mathcal K}$ there exists a fundamental domain $\Delta_{mod}$ for 
$PSL(2,{\mathbb Z})$ on ${\mathbb H}$ such that there is a conformal homeomorphism
$\varphi_a$ from $\Delta_{corr}=\Delta_{Cov} \cap \Delta_J$ to $\Delta_{mod}$ which 
(i) sends the vertices of the `croissant' 
$\Delta_{corr}$ to $0,\infty, i$ and $(-1+i\sqrt{3})/2$, and (ii) is equivariant 
with respect to the side-pairings induced by  $Cov$ and $J$ on the boundaries of 
$\Delta_{corr}$, and the side-pairings
$\rho:z \to -1-1/z$ and $\sigma:z \to -1/z$ on the boundaries of $\Delta_{mod}$. 
The existence of such a $\Delta_{mod}$ and (unique) $\varphi_a$ follows from the fact that 
the quotient orbifolds $\Delta_{corr}/{\mathcal F}_a$ and ${\mathbb H}/PSL(2,{\mathbb Z})$ are
(uniquely) conformally isomorphic, each being conformally a sphere with two cone points and a puncture. But note that the shape of the tile 
$\Delta_{mod}$ will vary with the
choice of Klein combination pair, and with $a$. 
We may extend the homeomorphism $\varphi_a$ equivariantly to a simply-connected union
 of `tiles' in $\Omega({\mathcal F}_a)$, including the one containing the critical value 
$v_a$ of ${\mathcal F}_a$, and  we can then define 
$$\Phi(a):=\varphi_a(v_a).$$
This map $\Phi$ can be shown to be well-defined on any simply-connected subset $U\setminus \mathcal{M}_{\Gamma}$
of ${\mathcal K}\setminus \mathcal{M}_{\Gamma}$ and to be a conformal homeomorphism from 
$U\setminus \mathcal{M}_{\Gamma}$ onto a
neighbourhood in ${\mathbb H}$ of the real interval $(-\infty,0)$ in $\partial {\mathbb H}$ 
(see \cite{BL4} for details). We remark that although a neighbourhood of 
$\mathcal{M}_{\Gamma}$ appears to be tiled in Figure \ref{mandelcorr}, the tile
boundaries here are not those of the pull-back via $\Phi$ of a tessellation of 
${\mathbb H}$ invariant under $PSL(2,{\mathbb Z})$: the figure is drawn
by plotting parameter points where ${\mathcal F}_a^n(v_a)$ lies on $\partial \Delta^{st}_{Cov}$ 
or $\partial \Delta^{st}_J$
for some $n\ge 0$, and these boundaries lift to different curves on ${\mathbb H}$ for different values of $a$. However the tile {\it vertices} in Figure \ref{mandelcorr} are the pull-backs of the {\it vertices} of a $PSL(2,{\mathbb Z})$-tessellation of ${\mathbb H}$.

\subsection{A Potential Function}\label{Green_fn}

In the Douady-Hubbard theory for quadratic polynomials a central role is played by a Green's function 
$G:{\mathbb C}\setminus K(Q_c) \to {\mathbb R}^{>0}$ for the filled Julia set,
which interacts with the dynamics via the formula $G(Q_c(z))=2G(z)$. We choose a potential function
for the limit set of ${\mathcal F}_a$ which will play an analogous role 
in our theory, although its interaction with 
the dynamics is more complicated (see Lemma \ref{alphabetapotential} below).\\

For every real $k>0$ the function $z \to k\log{|z|}$, from ${\mathbb C}\setminus {\mathbb D}$ to ${\mathbb R}^{\ge 0}$, 
is harmonic (since $\log{|z|}$ is the real part of $\log{z}$) and takes the value $0$ precisely on the unit circle $S^1$, the 
boundary of ${\mathbb D}$. The equipotentials are the circles $C_R$, with centre the origin and radius $R>1$.
We extend $z \to k\log{|z|}$ to a continuous function 
$\widehat \C\setminus {\mathbb D} \to {\mathbb R}^{\ge 0}\cup \{\infty\}$ by sending $\infty\in \widehat \C$ to 
$\infty\in \widehat \R$. \\

Let  $M: {\mathbb H} \to \widehat \C\setminus \overline{\mathbb D}$ denote the conformal homeomorphism $z\to \zeta=M(z)$, where
$$M(z)=\frac{z+i}{z-i}.$$

\begin{defn}\label{Green}
Let $\Psi:\widehat \C\setminus \overline{\mathbb D} \to {\mathbb R}^{> 0}\cup\{\infty\}$ be the function 
$\zeta \to (\log{|\zeta|})/2$, let
$\chi:{\mathbb H}\cup\{\infty\}\to {\mathbb R}^{> 0}\cup\{\infty\}$ be the composition $\chi=\Psi\circ M$, and 
for $a \in {\mathcal K}$ define $$G:\widehat \C\setminus \Lambda \to {\mathbb R}^{>0}\cup\{\infty\}$$
to be the composition $G=\chi\circ\varphi_a$, where $\varphi_a$ is the B\"ottcher map defined 
in Proposition \ref{Boettcher-prop}. 
\end{defn}

Note that:

(1) $G$ is harmonic on $\widehat \C\setminus \Lambda$; 

(2) setting $G(z)=0\ \forall z\in \Lambda$ extends $G$ to a continuous function $\widehat \C \to 
{\mathbb R}^{\ge 0}\cup\{\infty\}$. \\

We remark that while the circles $M^{-1}(C_R)$ in ${\mathbb H}$ 
are not the level sets of the `height' function $h(z)=y$ (for $z=x+iy\in{\mathbb H}$), they are close to these level sets
in the following sense. For points on the imaginary axis, we have
$$\chi(iy)=y+\frac{y^3}{3}+\frac{y^5}{5}+\frac{y^7}{7}\ldots$$
so for points $iy$ with $y$ small the value of $\chi$ is close to that of the height function (this is the reason for choosing $k$ to be $1/2$ in
the definition of $\Psi$). 
More generally
for points $z=x+iy \in {\mathbb H}$ which have $y$ small and $|x|$ in a bounded interval, $h(z)$ differs from $\chi(z)$ by a bounded factor. \\

We record the following relationship between $\chi(z)$, $\chi(\alpha z)$ and $\chi(\beta z)$ for points $z\in{\mathbb H}$ which are close to
the negative half of the real axis:

\begin{lemma}\label{alphabetapotential}
(i) There exist a constant $\lambda_+>1$ and neighbourhood $N\subset {\mathbb H}$ of the interval $[-\infty,0] \subset \partial{\mathbb H}$
such that $\chi(\alpha z)<\lambda_+\chi(z)$ and $\chi(\beta z)<\lambda_+\chi(z)$ for all $z \in N$. 

(ii) For every real $K>1$, there exist a constant $1<\lambda_- <\lambda_+$ and neighbourhoods $N_1\subset {\mathbb H}$ of the interval $[-K,-1] \subset \partial{\mathbb H}$ and  $N_2\subset {\mathbb H}$ of the interval $[-1,-1/K] \subset \partial{\mathbb H}$, such that 
$\chi(\alpha z)>\lambda_-\chi(z)$ $\forall z \in N_1$, and $\chi(\beta z)>\lambda_-\chi(z)$ $\forall z\in N_2$.

\end{lemma}

{\bf Proof}

Conjugation by $M:z \to \zeta=(z+i)/(z-i)$ sends the automorphism $\alpha: z \to z+1$ of ${\mathbb H}$ to the automorphism 
$\zeta \to M\alpha M^{-1}(\zeta)$ of 
$\widehat \C\setminus \overline{\mathbb D}$. \\

Since $\chi(z)=(\log{|\zeta|})/2$, where $\zeta=M(z)$, 
to prove (i) we must show that there exists a bound $\lambda_+$ such that every point $\zeta$ in the complement of the unit 
disc which is sufficiently close to $M[-\infty,0]$ has 
$$\log(|M\alpha M^{-1}(\zeta)| )< \lambda_+\log(|\zeta|)\ \ \ {\rm and} \ \ \ \log(|M\beta M^{-1}(\zeta)|) < \lambda_+\log(|\zeta|).$$
However $M[-\infty,0]$ consists of the lower half of the unit circle, traversed clockwise, and so the 
points of $\widehat \C\setminus \overline{\mathbb D}$
close to $M[-\infty,0]$ are the points $\zeta=Re^{-i\theta}$ which have $R=1+r$ with $r>0$ small, and also have $0 \le \theta \le \pi$. 
But, if we neglect $r^2$, every point $\zeta$ with $|\zeta|=R=1+r$ has $\log(|\zeta|)\sim r=|\zeta|-1$, and as the images of $\zeta$
under $M\alpha M^{-1}$ and $M\beta M^{-1}$ are also close to the unit circle we are reduced 
to proving that there exists a constant $\lambda_+>0$ such that
for $\zeta=Re^{-i\theta}$, with $R=1+r$ and $0 \le \theta \le \pi$, we have
$$(|M\alpha M^{-1}(\zeta)|-1)< \lambda_+ r\ \ \ {\rm and} \ \ \ (|M\beta M^{-1}(\zeta)|-1) < \lambda_+r.$$
However
$$M\alpha M^{-1}(Re^{-i\theta})=\frac{-(1+2i)R^{-i\theta}+1}{-Re^{-i\theta}+(1-2i)},$$
so
$$|M\alpha M^{-1}(Re^{-i\theta})|^2=\frac{5R^2-2R(\cos{\theta} +2\sin{\theta})+1}{5-2R(\cos{\theta}+2\sin{\theta})+R^2},$$
which, setting $R=r+1$ and assuming $r^2$ to be negligible, simplifies to
$$ \frac{6-2(\cos{\theta}+2\sin{\theta})+r(10-2(\cos{\theta}+2\sin{\theta}))}{6-2(\cos{\theta}+2\sin{\theta})+r(2-2(\cos{\theta}+2\sin{\theta}))}
\sim 1+ \frac{4r}{3-(\cos{\theta}+2\sin{\theta})}.$$
Thus (still assuming $r^2$ negligible)
$$|M\alpha M^{-1}(Re^{-i\theta})|-1=\frac{2r}{3-(\cos{\theta}+2\sin{\theta})}.$$
When $\tan{\theta}=2$ this function of $\theta$ attains its maximum value
$$\frac{2r}{3-\sqrt{5}}.$$
It follows that $M\alpha M^{-1}$ (and hence also $\alpha$) 
increases potential by a factor of at most $$\frac{2}{3-\sqrt{5}}=(3+\sqrt{5})/2.$$ A similar computation gives
the same upper bound for the multiplier of $M\beta M^{-1}$ (and hence also $\beta$) on potential. In both computations we  
assume $r^2$ to be
negligible, so they are only valid in the limit as we approach the unit circle. But we may obtain a 
bound $\lambda_+$ which 
is valid on a sufficiently small neighbourhood of the unit circle by setting
$$\lambda_+=((3+\sqrt{5})/2)+\varepsilon$$
for any $\varepsilon>0$. 
\\

To prove the statement in (ii) concerning $\alpha$ we must find a constant $\lambda_->1$ such that for every point 
$\zeta$ in the complement of the unit disc which is sufficiently close to $M[-K,-1]$ we have
$$(|M\alpha M^{-1}(\zeta)| -1)> \lambda_-(|\zeta|-1)$$
For this we need a {\it lower} bound  ($>1$) for
$$\frac{2r}{3-(\cos{\theta}+2\sin{\theta})}.$$
However for $0<\theta\le \pi/2$ we have $$\cos{\theta}+2\sin{\theta}>1,$$ and so provided we bound $\theta$ 
away from $0$, so that say $0<\delta<\theta<\pi/2$ for some small constant $\delta$, we can find a constant $k$ such that
 $$\cos{\theta}+2\sin{\theta}>k>1$$
 for all $\theta\in (\delta,\pi/2)$, and thus we can find a constant $\lambda_->1$ with the desired property.
 Finally, the part of (ii) concerning $\beta$ follows from a calculation for $M\beta M^{-1}$, which for $\pi/2 \le \theta <\pi$ 
 gives the same constant.
 \qed

\section{Geodesics}\label{geod}
When $a\in {\mathcal C}_\Gamma$, so that $\Lambda({\mathcal F}_a)$ is connected, we can use the B\"ottcher isomorphism $\varphi_a$ from 
$\Omega({\mathcal F}_a)=\widehat \C\setminus \Lambda({\mathcal F}_a)$ to ${\mathbb H}$
to pull back the hyperbolic metric on ${\mathbb H}$ to the hyperbolic metric on $\Omega({\mathcal F}_a)$.
The geodesics for this metric become the analogues in our setting of the {\it external rays}
 defined by Douady and Hubbard
for quadratic polynomials with connected Julia sets.
In \ref{per-geo} we associate periodic geodesics to (finite) words $W$, and we prove that these geodesics 
land at both ends (Proposition \ref{periodic_rays_land}). Generalising, in \ref{bi-inf} we associate geodesics to {\it bi-infinite sequences} $S\in \{\alpha,\beta\}^{\mathbb Z}$, 
and in \ref{Mink_qm} we explain how this association is a manifestation of Minkowski's {\it question mark map}. 
In \ref{St} we consider 
{\it Sturmian sequences}: thought of as binary representations of real numbers these have  orbit under the doubling map 
arranged in the same order as a rigid rotation of the circle ${\mathbb R}/{\mathbb Z}$.
We conclude the section by computing { \it bounds on the multipliers} of Sturmian words and sequences.

\subsection{Geodesics associated to finite words}\label{per-geo}
As before we write $\alpha$ for $z\to z+1$ and $\beta$ for $z\to z/(z+1)$, now both acting on the boundary 
$\widehat \R={\mathbb R}\cup\{\infty\}$ of ${\mathbb H}$ as well as on ${\mathbb H}$ itself. We recall that $\beta$ is conjugate
to $\alpha^{-1}$ (via the inversion $z \to -1/z$) and that both $\alpha$ and $\beta$ are parabolic, having unique fixed point $\infty$ and $0$ respectively. A (finite) word 
$W$ in the letters $\alpha$ and $\beta$ acts on the left on points $z\in {\mathbb H}$ (for example $W=\alpha^2\beta^2\alpha$
acts by $z \to \alpha^2(\beta^2(\alpha z))$).

\begin{lemma} \label{per_geo}

(i) For every word $W$ consisting of a finite sequence of positive powers of both $\alpha$ and $\beta$, 
there is a unique $W$-invariant geodesic $\gamma(W)$ in $\mathbb H$.

(ii) The end points of $\gamma(W)$ are a repelling fixed point $x^-(W)$ in ${\mathbb R}^{<0}$ and an attracting
fixed point $x^+(W)$ in ${\mathbb R}^{>0}$.

\end{lemma}

{\bf Proof} 

The equation $Wz=z$ has the form $(az+b)/(cz+d)=z$, with $a,b,c,d$ non-negative integers, so it 
has two real solutions if $W$ is hyperbolic, one real solution if $W$ is parabolic and no real solution if $W$ is elliptic. 
We first note that the hypothesis that the word $W$ contains positive powers of {\it both} $\alpha$ and $\beta$ is necessary.
The words $\alpha^n$ and $\beta^n$ ($n>0$), being parabolic, have just one fixed point each and do not fix any geodesic in $\H$.\\

(i) The products $\alpha\beta$ and $\beta\alpha$ both have trace $3$. Building up $W$ inductively from either $\alpha\beta$ or
$\beta\alpha$ by multiplying on one side or the other by $\alpha$ or $\beta$ either leaves the trace unchanged or increases 
it, since the matrices being multiplied together have all of their entries non-negative.  We
deduce that $W$ is hyperbolic, and so it  has two fixed points, both real.
As the product $-b/c$ of these two real numbers is negative, they lie on either side on $0$. 
Denote these fixed points by $x^-(W)$ and $x^+(W)$ and let $\gamma(W)$ denote the geodesic in 
${\mathbb H}$ which joins them, the {\it axis} of the loxodromic M\"obius transformation $W$.\\

(ii) It will suffice to show that the derivative of $W$ at its fixed point $x^+(W)$ has modulus less than one, i.e. 
that this point is the attractor for the action of $W$ on $\gamma(W)$. But this derivative is $1/|cx^+(W)+d|^2$, 
which can also be written as $(|x^+(W)|/|ax^+(W)+b|)^2$ (since $Wx^+(W)=x^+(W)$) and this is clearly less than 
$1$ as $a$ and $b$ are positive integers.
$\square$\\

\begin{example}\label{AAB}
$W=\alpha\alpha\beta$. 

$$\alpha\alpha\beta(z)=\frac{z}{z+1}+2=\frac{3z+2}{z+1}$$
so $x^-(W)=-(\sqrt{3}-1)$ and $x^+(W)=\sqrt{3}+1$.

The orbit of $x^-(W)$ is the cycle $$x^-(W) \to  \beta x^-(W) \to \alpha\beta x^-(W) \to \alpha\alpha\beta x^-(W)=x^-(W),$$ 
that is, $P_0=-(\sqrt{3}-1) \to P_1=-(\sqrt{3}+1) \to P_2=-\sqrt{3}\to P_0$.

The orbit of $x^+(W)$ is the cycle $$x^+(W) \to  \beta x^+(W) \to \alpha\beta x^+(W) \to \alpha\alpha\beta x^+(W)=x^+(W),$$ 
that is, $Q_0=\sqrt{3}+1\to Q_1=\sqrt{3}-1 \to Q_2=\sqrt{3} \to Q_0$. 

\end{example}

\begin{figure}
\begin{center}
\scalebox{1.0}{\includegraphics{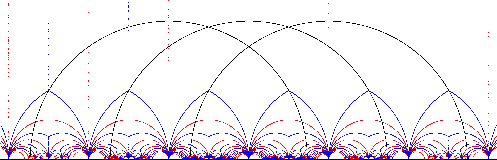}}
\caption{\small The geodesics in ${\mathbb H}$ in Example 1 (the red vertical lines are 
at integer values on the real axis, from $-3$ to $+3$).}\label{period3example}
\end{center}
\end{figure}

See Figure \ref{period3example} for an illustration of these geodesics in ${\mathbb H}$. \\

\begin{prop}\label{periodic_rays_land}
For every $a\in {\mathcal C}_\Gamma$ and every word $W$ consisting of a finite sequence of positive powers of 
both $\alpha$ and $\beta$, the inverse  of the B\"ottcher map
$$\psi_a:{\mathbb H} \to \Omega$$
extends continuously to $x^-(W)$ and $x^+(W)$.
Under this extension,
\begin{itemize}
 \item $\psi_a(x^-(W)) \in \partial \Lambda_-$ is a repelling 
periodic point of the $2$-to-$1$ restriction $f_a$ of ${\mathcal F}_a$ defined on a neghbourhood of $\Lambda_-$;
\item $\psi_a(x^+(W)) \in \partial \Lambda_+ $ is a repelling 
periodic point of the restriction $Jf_aJ$ of
${\mathcal F}_a^{-1}$ defined on a neighbourhood of $\Lambda_+$, and thus an attracting periodic point of ${\mathcal F}_a$.
\end{itemize}
\end{prop}

{\bf Proof}

We follow the method of proof of Theorem 18.10 in \cite{M}. First we modify our earlier definition of the potential function 
$G$ (Definition \ref{Green}, Section \ref{Green_fn}), by precomposing the map $\chi$ in that definition by a M\"obius transformation 
which sends the geodesic $\gamma(W)$ 
to the imaginary axis, with the initial end of $\gamma(W)$ going to $0$ and the final end to $\infty$. 
Denote this new potential function by $G_W$.\\

Now divide the imaginary axis in ${\mathbb H}$ into segments of Poincar\'e
length $\ln(\mu(W))$, where $\mu(W)$ is the multiplier of $W$ (the square of the eigenvalue which is  $>1$).
Correspondingly, parametrise $\varphi_a^{-1}(\gamma(W))$ as a path $p:{\mathbb R}  \to \widehat \C\setminus \Lambda$ 
with ${\mathbb R}$ divided into unit intervals $I_k$ ($k\in{\mathbb Z}$) each mapped isometrically 
to the next by $f_a^q$, where $q$ is the length of the word $W$.
Since $G_W(p(s))$ tends to zero as $s\to -\infty$ any limit point $\hat{z}$ of $\{p(s):s\le 0\}$ must belong to 
$\partial \Lambda_-$. Following the same reasoning as in the proof of Theorem 18.10 in \cite{M}
we deduce that the geodesic $\varphi_a^{-1}(\gamma(W))$ has a limit point $\hat{z}$ at its initial end,
 that this limit point is a fixed point 
of $f_a^q$, that the geodesic lands at $\hat{z}$, and that $\hat{z}$ is necessarily repelling or parabolic.\\

Conjugating ${\mathcal F}_a$ by $J$ sends ${\mathcal F}_a$ to ${\mathcal F}_a^{-1}$ and interchanges $\Lambda_-$
with $\Lambda_+$. The result for $x^+(W)$ follows.
$\square$

\begin{remark}
As $J\alpha J=\beta^{-1}$ and $J\beta J=\alpha^{-1}$, the orbit of $\psi_a(x^+(W))$ on 
$\Lambda_+({\mathcal F}_a)$ is not the
$J$-image of the orbit of $\psi_a(x^-(W))$ on $\Lambda_-({\mathcal F}_a)$ for most words 
$W$ (an exception being $W=\alpha\beta$). 
As we vary the parameter $a$, the points of a periodic orbit may collapse together.
For example, consider $a=4.53926-0.439437i$ 
corresponding to a mating between $PSL(2,{\mathbb Z})$ and the (superattracting) `co-rabbit'. 
Here the orbit 
$\{P_0,P_1,P_2\}$ of $x^-(\alpha\alpha\beta)$ 
(the ends of geodesics to the left of the origin in Figure \ref{period3example}) are identified 
under $\psi_a$ to a single point, but the points of the orbit of $x^+(\alpha\alpha\beta)$ are not:
at this value of $a$ the three periodic geodesics emanating from the fixed point in 
$\Lambda_-({\mathcal F}_a)$ have their opposite ends at distinct points of 
$\Lambda_+({\mathcal F}_a)$
(the points of the `rabbit orbit', not the `co-rabbit orbit').

\end{remark}

\begin{cor}
For $a\in \mathcal{M}_{\Gamma}$ the inverse $\psi_a:{\mathbb H} \to \Omega$ of the B\"ottcher map $\varphi_a$
extends continuously to all quadratic irrationals in ${\mathbb R}$.
\end{cor}

{\bf Proof}

Recall that a {\it quadratic irrational} is a root of a quadratic equation which has integer coefficients and real 
but not rational solutions, and that the positive quadratic irrationals are precisely the positive
real numbers which have periodic or pre-periodic continued fraction expansions. We already know
that periodic rays land, so it just remains to consider the strictly preperiodic case. Given any positive
quadratic irrational $x$, the periodic `tail' of the continued fraction is an $x^+(W)$ for some $W$, and so $x$ 
can be written $W'x^+(W)$ for some finite word $W'$ in $\alpha$ and $\beta$. The map $\psi_a$ extends 
continuously to $x^+(W)$ by Proposition \ref{periodic_rays_land},
and so it extends to every point on the orbit of $x^+(W)$ under the correspondence, in particular to $W'x^+(W)$.
As $x^-(W)=Jx^+(W'')J$ for some $W''$ (usually different from $W$) the statement for negative quadratic
irrationals follows. $\square$

\subsection{Geodesics associated to bi-infinite sequences}\label{bi-inf}
Let $S\in \{\alpha,\beta\}^{\mathbb Z}$. We think of $S$ as 
a bi-infinite sequence
$$\ldots g_{n} \ldots g_1g_0\ \cdot \ g_{-1}\ldots g_{-n}\ldots$$
where each $g_i$ is $\alpha$ or $\beta$, and the dot between $g_0$ and $g_{-1}$ is a position marker. 
In this section we will associate to $S$ a geodesic 
$\gamma(S)$, having left hand end-point in ${\mathbb R}^{\le 0}$ and right hand end-point in ${\mathbb R}^{\ge 0}$, and 
we shall see (Lemmas \ref{W} and \ref{shift} below) that this process generalises
our earlier definition for a finite word.\\

Collecting up contiguous occurrences of each letter, we can write $S$ in the form:
$$S=\ldots    \beta^{m_3}\alpha^{m_2}\beta^{m_1}\alpha^{m_0}\ \cdot\ \alpha^{n_0}\beta^{n_1}\alpha^{n_2}\beta^{n_3}\ldots$$
where $m_0,n_0 \in \{0\}\cup {\mathbb N} \cup \{\infty\}$and  $m_i,n_i \in {\mathbb N} \cup \{\infty\}$ for $i>0$ (and 
if any $m_i$ or $n_i$ is $\infty$ then subsequent $m_j$  or $n_j$ are undefined).\\

We may regard $S$ as a back-to-back pair of symbol sequences defining
a pair of points $x^-(S)\in {\mathbb R}^{\le 0}$ and $x^+(S)\in {\mathbb R}^{\ge 0}$, and thereby associate to $S$ the \textit{directed geodesic} $\gamma(S)$, which has ends as follows: 
$$x^-(S)= -\infty \mbox{ when } m_0=\infty, \mbox{ and } x^-(S)=0 \mbox{ when } m_0=0 \mbox{ and } m_1=\infty$$
$$x^-(S)=-[m_0;m_1,m_2,m_3,\ldots] \in {\mathbb R}^{<0}, \mbox{ in all the other cases};$$
$$x^+(S)= \infty \mbox{ when } n_0=\infty, \mbox{ and } x^+(S)=0 \mbox{ when } n_0=0 \mbox{ and } n_1=\infty$$
$$x^+(S)=[n_0;n_1,n_2,n_3,\ldots] \in {\mathbb R}^{>0}, \mbox{ in all the other cases}.$$

Note that there are 
dynamical systems associated to $S$
which have these points as attractors. Starting from the observation that $\alpha^{-1}({\mathbb R}^{<0})=(-\infty,-1)$ 
and $\beta^{-1}({\mathbb R}^{<0})=(-1,0)$ are disjoint intervals in ${\mathbb R}^{<0}$ whose closures cover 
${\mathbb R}^{<0}$, we can find an infinite composition of applications of $\alpha^{-1}$ and $\beta^{-1}$ for which 
the images of ${\mathbb R}^{<0}$ nest down to any chosen point of ${\mathbb R}^{<0}$ (this is what 
we do when we write down a continued fraction expression for the chosen point). Similarly any chosen 
point of ${\mathbb R}^{>0}$ is the limit of the images of a sequence of compositions of applications
of $\alpha$ and $\beta$ to ${\mathbb R}^{>0}$. Formally, for $n> 0$ let $G_n$ denote 
the composite M\"obius transformation $g_0^{-1}g_1^{-1}\ldots g_n^{-1}$, 
and $G_{-n}$ denote $g_{-1}g_{-2}\ldots g_{-n}$.
Then for any choice of $z \in {\mathbb H}\cup{\mathbb R}^{<0}$ we have:
$$ \lim_{n\to \infty} G_n(z)  =x^-(S) \in {\mathbb R}^{\le 0}\cup\{-\infty\}$$
and for any choice of $z \in {\mathbb H}\cup{\mathbb R}^{>0}$ we have:
$$\lim_{n\to \infty} G_{-n}(z)  =x^+(S) \in {\mathbb R}^{\ge 0}\cup\{\infty\}.$$

\begin{lemma}\label{W}
If $W$ is a finite word in $\alpha$ and $\beta$, and $S(W)=\overline{W}\ \cdot \ \overline{W}$ (the bi-infinite sequence
consisting of repeats of $W$), then
$$\gamma(S(W))=\gamma(W)$$
where $\gamma(W)$ is as in Lemma \ref{per_geo}. 
\end{lemma}

{\bf Proof}

Immediate from definitions. $\square$

\begin{lemma}\label{shift} Let $S=(g_n)_{n\in {\mathbb Z}}$, and 
let $\sigma$ denote the operation of moving 
the position marker in $S$ one place to the left (so $(\sigma S)_n=g_{n-1}$). Then
$$\gamma(\sigma(S))=g_0(\gamma(S)).$$
\end{lemma}

{\bf Proof}

This follows at once from the algorithms for $x^+(S)$ and $x^-(S)$ above, together with the observation 
that $x\to -x$ conjugates $\alpha$ and $\beta$ on ${\mathbb R}^{\ge 0}$ to $\alpha^{-1}$ and $\beta^{-1}$ respectively on ${\mathbb R}^{\le 0}$. $\square$\\

{\bf Example \ref{AAB} revisited.} The bi-infinite sequence corresponding to the finite word $W=\alpha\alpha\beta$ of Example \ref{AAB} is:
$$S=  \ldots \alpha^2\beta     \alpha^2\beta\ \cdot\ \alpha^2\beta \alpha^2\beta \ldots$$
so the definition above gives:
$$x^-(S)=-[0;1,2,1,\ldots]=-(\sqrt{3}-1); \quad x^+(S)=[2;1,2,1,\ldots]=\sqrt{3}+1$$
which agrees with our calculation in Example \ref{AAB}. Furthermore
$$\sigma(S)= \ldots \alpha^2\beta \alpha^2\beta\alpha^2\ \cdot\ \beta \alpha^2\beta \alpha^2 \beta \dots$$
which gives
$$x^-(\sigma(S))=-[2;1,2,1,2,\ldots]=-(\sqrt{3}+1);\quad x^+(\sigma(S))=[0;1,2,1,2,\ldots]=\sqrt{3}-1$$
confirming that $\gamma(\sigma(S))$ is indeed equal to $\beta\gamma(S)$ in this example.

\subsection{Minkowski's `question mark' map}\label{Mink_qm}

In the preceding subsection we associated continued fraction expressions representing real numbers $x^-(S)$ and 
$x^+(S)$ to each marked bi-infinite sequence $S$ of `$\alpha$'s and `$\beta$'s. Equally, we may associate binary expressions 
representing real numbers in the interval $[0,1]$ to the sequences $n_0,n_1,n_2,\ldots$ and $m_0,m_1,m_2,\ldots$ which code $S$, 
now with the symbols `$1$' and `$0$' corresponding to `$\alpha$' and `$\beta$' respectively. Indeed the correspondence between 
real numbers expressed as continued fractions, and real numbers in the interval $[0,1]$ expressed in binary, is at the heart of
the existence of matings between $PSL(2,{\mathbb Z})$ and quadratic polynomials in \cite{BP}. 
The key is Minkowski's `question mark' map  (see \cite{Mi} p171-172).\\

Minkowski's map is a homeomorphism from the unit interval $(0,1)\subset {\mathbb R}$
to itself. We consider the following slightly modified version (which we denote by `?' as did Minkowski his map). Our map
is the homeomorphism from 
$(0,\infty)={\mathbb R}^{>0}$ to the unit interval:
$$?([a_0;a_1,a_2\ldots])=0.\underbrace{1\ldots 1}_{a_0}\underbrace{0\ldots 0}_{a_1}\underbrace{1\ldots 1}_{a_2}\ldots$$
where on the left hand side is a continued fraction expression and on the right-hand side is a binary expression.
Given any point $x^{-}(S)\in {\mathbb R}^{<0}$ the point $?(-x^{-}(S))\in [0,1]$ is the real number which has the binary 
expression $0.t_0t_1t_2\ldots$
where each $t_j$, $j\ge 0$ is $0$ or $1$ according to whether $g_j$ is $\beta$ or $\alpha$ in 
the bi-infinite word $S$. The bijection
$$ 0.t_0t_1t_2\ldots   \leftrightarrow   0.g_0g_1g_2\ldots\  \ \mbox{ given by } \ 1 \leftrightarrow \alpha,\ 0\leftrightarrow \beta$$
gives us an explicit formula for the correspondence between periodic external rays for a quadratic polynomial $Q_c$ and periodic
geodesics for a mating between $Q_c$ and $PSL(2,{\mathbb Z})$. For example, the external ray labelled $.110110110\ldots$ and its
orbit under the doubling map (which has rotation number $2/3$), correspond to the periodic geodesic 
$\gamma(W)$, $W=\alpha\alpha\beta$, of 
Example \ref{AAB}, and its orbit under the cycle `apply $\beta$ then $\alpha$ then $\alpha$':
$$3/7=.011011\ldots \leftrightarrow [0;1,2,1,2,\ldots]=\sqrt{3}-1=-x^-(\alpha\alpha\beta);$$
$$6/7=.110110\ldots \leftrightarrow [2;1,2,1,2,\ldots]=\sqrt{3}+1=-x^-(\beta\alpha\alpha);$$
$$5/7=.101101\ldots \leftrightarrow [1;1,2,1,2,\dots]=\quad \sqrt{3}\quad  =-x^-(\alpha\beta\alpha).$$

\subsection{Sturmian sequences, rotation numbers, bounds on multipliers}\label{St}
There are many equivalent definitions of the term `Sturmian' (which is due to Morse and Hedlund); it is usually applied to infinite
sequences, but can also be applied to bi-infinite sequences or to (finite) words.
For the purposes of the current article:

\begin{defn}
An infinite or bi-infinite sequence in the symbols $0$ and $1$ is {\it Sturmian} if for each $n\in {\mathbb N}$ the numbers of
$1$'s in any two blocks of $n$ consecutive symbols differ by at most $1$. A word $W$ in $0$'s and $1$'s is {\it Sturmian}
if the infinite sequence $\cdot \overline{W}$ is Sturmian (or equivalently the bi-infinite sequence 
$\overline{W} \cdot \overline{W}$) is Sturmian. Here, as usual, $\overline{W}$ denotes a repeated sequence of blocks $W$.
\end{defn}

Sturmian sequences are also known as `balanced sequences'. 
They occur as the `cutting sequences' of straight lines of rational or 
irrational slope on an integer grid of squares (with appropriate conventions where lines pass through vertices: see the 
proof of the Proposition below). A useful way to 
characterise an infinite Sturmian sequence is as a sequence of $0$'s and $1$'s such that the real number in $[0,1]$ it
represents in binary has orbit under the doubling map a sequence of points arranged in the same order around the 
circle ${\mathbb R}/{\mathbb Z}$ as for a rigid rotation (see \cite{BS}). A bi-infinite word is Sturmian if and only if
the infinite sequence to the right (equivalently left) of the marker has this property wherever the 
marker is placed. This allows us to definite a {\it rotation number} for such words.

\begin{defn}\label{rotst}
 The rotation number of a Sturmian sequence is that of the corresponding rigid rotation of the circle: equivalently it is 
the limiting frequency with which the digit
`$1$' appears in subwords of length $n$, as $n$ tends to $\infty$. 
\end{defn}

The following Proposition lists the properties of Sturmian sequences that we shall make use of in our proofs in subsequent Sections: 

\begin{prop}\label{Sturm}

\begin{enumerate}
\item[(i)]
For each rational $0<p/q<1$ there is exactly one Sturmian word $T_{p/q}$ (up to cyclic equivalence) of 
length $q$ and rotation number $p/q$, and there are two bi-infinite non-periodic Sturmian sequences of rotation 
number $p/q$. Both ends of each of the non-periodic Sturmian sequences consist of repeated copies of $T_{p/q}$.

\item[(ii)]
For each irrational $0<\nu<1$, the set of Sturmian sequences of rotation number $\nu$ forms a 
Cantor set $C_{\nu}$ contained in the real interval $(0,1)$, when these sequences are regarded as binary expressions for real numbers. 
For each bi-infinite Sturmian sequence $S$ of rotation number $\nu$, and each shift (left or right), the orbit of $x^-(S)$ 
under the shift is dense in $C_\nu$, as is that of $x^+(S)$.

\end{enumerate}
\end{prop}

{\bf Proof}

We omit details, but these properties follow from the characterisation of a Sturmian sequence of rotation 
number $0 \le \nu <1$ as a sequence obtained by 
the `staircase algorithm' \cite{BS}, applied to a straight line $L$ of slope $\nu$ superimposed on an integer grid of lines.
This algorithm codes the maximum integer staircase that fits below $L$ by writing
`$0$' for a horizontal move, and `$1$' for a horizontal plus vertical move. If $L$ passes through 
one or more vertices the rule to obtain a staircase is to cut $L$ at a point where it does not meet 
a grid line and then parallel translate the two halves of $L$ infinitesimally in opposite directions so that they 
no longer pass through any vertices. Thus for each irrational $\nu$ there is a countable set of Sturmian 
sequences of rotation number $\nu$ which have two continuations to bi-infinite Sturmian sequences (these correspond
 to the lines $L$ of slope $\nu$ which have continuations to the left which pass through a vertex), whereas 
all other Sturmian sequences of rotation number $\nu$ have unique continuations to bi-infinite Sturmian 
sequences. The result concerning density of every orbit in $C_{\nu}$, under either
shift, follows from the fact that the intersections between $L$ and vertical grid lines, when projected to the vertical
axis ${\mathbb R}$ and then to the circle ${\mathbb R}/{\mathbb Z}$, become the points of an orbit of an irrational
rigid rotation of the circle, and thus dense in the circle under either forward or backward  iteration.
\qed

\begin{example}
(a) $T_{1/3}=001$. The two non-periodic bi-infinite Sturmian sequences of rotation number $1/3$ are
$\overline{(001)}(01)\overline{(001)}$ and $\overline{(001)}(0001)\overline{(001)}$.

(b) The limit of the Sturmian words $10$, $101$, $10110$,  $10110101$, $\dots\ $ (obtained from $10$ by repeatedly applying
the substitutions $1 \to 10$, $0\to 1$) is a Sturmian sequence of 
rotation number the golden mean. This has two continuations to the left
yielding bi-infinite Sturmian sequences, namely $\dots 10110101\cdot 10110101\ldots$ and $\dots 10110110\cdot 10110101\ldots$

\end{example}

\begin{cor}\label{invariant_sets_finite} If $A\subset{\mathbb R}/{\mathbb Z}$ is a closed invariant Sturmian subset 
on which the doubling map acts injectively, then the rotation number $\nu$ of the doubling map restricted to $A$ is a rational $p/q$, and the points of $A$ 
are the real numbers whose binary expressions are the cyclic permutations of $T_{p/q}$. In particular $A$ is finite.
\end{cor}

{\bf Proof}

Let $x\in A$. Since the doubling map restricted to $A$ is a homeomorphism we can continue 
the binary expression for $x$ to a bi-infinite Sturmian sequence.
First suppose that $\nu$ is rational.  If the bi-infinite sequence for $x$ is either of the two non-periodic sequences listed 
in part (i) of the Proposition, the orbit of $x$ under the doubling map is not injective, since it contains a point outside the orbit of 
$(T_{p/q})^\infty$ mapping onto this orbit.
Now suppose that $\nu$ is irrational. In this case by part (ii) of the Proposition, since $A$ is closed we deduce that $A$ is
the Cantor set $C_\nu$. But the doubling map on $C_\nu$ is non-injective, since it sends the two ends of the longest gap in $C_\nu$ 
to a single point (the longest gap has length $1/2$, see \cite{BS}). \qed\\

We now replace `$0$' by `$\beta$' and `$1$' by `$\alpha$', and consider the finite Sturmian block $T_{p/q}$ as a composition of 
$q$ matrices, $p$ of which are copies
of $\alpha$ and $q-p$ of which are copies of $\beta$. 
We shall establish upper and lower bounds on the multiplier $\mu_{p/q}$ of $T_{p/q}$ at $x^-(T_{p/q})$, as a consequence of the 
following result:

\begin{prop}\label{multipliers} Let $r>1$. 
\begin{enumerate}
\item[(i)]
If $W=\alpha^{r-1}\beta$ then the multiplier $\mu(W)$ of $W$ satisfies the inequality:
$$r^{2}<\mu(W)<(r+1)^{2}.$$
\item[(ii)]
If $W$ is a word made up of $s>1$ blocks, each of the form either $\alpha^{r-1}\beta$ or $\alpha^r\beta$, then 
the multiplier $\mu(W)$ of $W$ satisfies the inequality:
$$r^{2s} < \mu(W) < (r+2)^{2s}.$$
\item[(iii)] The inequalities above also hold when $\alpha^{r-1}\beta$ and $\alpha^r\beta$ are replaced by $\beta^{r-1}\alpha$ and $\beta^r\alpha$.
\end{enumerate}
\end{prop}

{\bf Proof}

We shall estimate the positions of points on the orbit  of $x^-(W)$ on ${\mathbb R}^{<0}$, and the value of the derivative
of $\alpha$ or $\beta$ (as appropriate) at each point, then multiply these derivatives together to get an estimate of the
multiplier of the orbit. As the derivative of $\alpha:x\to x+1$ is $1$ everywhere on ${\mathbb R}$, we only have to compute 
the derivatives at orbit points where the map being applied is $\beta:x \to x/(x+1)$.\\

(i) The unique point of the orbit of $x^-(W)$ which lies in the open interval $(-1,0)$
is the solution $x_r<0$
of the quadratic equation $\alpha^{r-1}\beta(x_r)=x_r$. We could solve this equation to find $x_r$, and then 
compute the derivative of $\beta$ at $x_r$ to find the multiplier of the orbit, but equally we can proceed as follows.
The trace of the matrix $\alpha^{r-1}\beta$ is $r+1$, so its eigenvalues are $\lambda$ and $\lambda^{-1}$ where
$$\lambda=\frac{(r+1)+\sqrt{r^2+2r-3}}{2}.$$
Thus there is a M\"obius conjugacy from $\alpha^{r-1}\beta$ to the map 
$z \to {\lambda z}/{\lambda^{-1}}=\lambda^2 z$, and
the derivative of $\alpha^{r-1}\beta$ at its expanding fixed point $x_r$ is therefore
$$\mu(W)=\lambda^2=\left(\frac{(r+1)+\sqrt{r^2+2r-3}}{2}\right)^2.$$
In particular 
$$r^2<\mu(W)<(r+1)^2.$$\\

(ii) The $s$ points of the orbit of $x^-(W)$ which lie in $(-1,0)$ all lie between $x_{r}$ and $x_{r+1}$, since 
the Minkowski question mark map preserves order and we know that a binary sequence made up of
blocks of the form `$1$ followed by $r-1$ copies of $0$', or `$1$ followed by $r$ copies of $0$', represents 
a real number which lies between the numbers represented by two periodic sequences made up of  
copies of just one of these blocks.
Finally, since the derivative of 
$\beta$ at $x$ is $1/(1+x)^2$, which is monotonic in $x$,
and the multiplier of the orbit is the product of the values of the derivative of $\beta$ at the points of the orbit
which are in $(-1,0)$, the general result follows from our initial special case calculation. \\

(iii) follows at once from the facts that $z\to -1/z$ conjugates $\alpha$ to $\beta^{-1}$ and $\beta$ to $\alpha^{-1}$ 
on ${\mathbb H}$, and that for any invertible matrix $M$ the eigenvalues of $M^{-1}$ are the inverses of
those of $M$.
\qed

\begin{cor}\label{estimates}  
\begin{enumerate}

\item[(i)]
For all $0< p/q \le1/2$, the multiplier $\mu_{p/q}$ of $T_{p/q}$ satisfies
$$\lfloor q/p \rfloor^{2p}<\mu_{p/q}<(1+\lceil q/p \rceil)^{2p}.$$

\item[(ii)]
For all $1/2\le p/q <1$, the multiplier $\mu_{p/q}$ satisfies
$$\lfloor q/(q-p) \rfloor^{2(q-p)}<\mu_{p/q}<(1+\lceil q/(q-p) \rceil)^{2(q-p)}.$$

\end{enumerate}

\end{cor}

{\bf Proof}
First suppose that $0<p/q\le 1/2$. The word $T_{1/q}$ is (up to cyclic equivalence) $\alpha^{q-1}\beta$, and, 
when $p>1$, $T_{p/q}$ is made up of blocks $\alpha^{r-1}\beta$ and $\alpha^r\beta$, where $r=\lfloor q/p \rfloor$
and $r+1=\lceil q/p \rceil$. \\

Now (i) and (ii) follow from  parts (ii) and (iii) of Proposition \ref{multipliers} respectively. 
\qed

\section{The proof of Theorem \ref{repelling_fp}}\label{repellor}
Motivated by results for polynomial maps
and having proved in Proposition \ref{periodic_rays_land} (Section \ref{per-geo}) that every periodic geodesic lands, we now set out to 
prove Theorem \ref{repelling_fp}, that for 
$f_a$ with $a \in {\mathcal C}_\Gamma$, every repelling periodic point in $\Lambda_-$ 
is the landing point of a periodic geodesic.\\

We start by establishing notation for linearisation around a repelling fixed point, and some 
preparatory results concerning its properties. 
Let $\hat{z} \in \Lambda_-$ be a repelling fixed point of $f$($=f_a$), and let $\omega$ denote the 
derivative of $f$ at $\hat{z}$  (so $|\omega|>1$). 
Let ${\mathbb D}$ be the open unit disc and $\lambda$ be a Koenigs linearisation, that is a 
conformal homeomorphism from ${\mathbb D}$ to an open topological disc containing $\hat{z}$, conjugating 
the map $\times \omega$ (multiplication by $\omega$) to $f$. 
Let $C_{-n}$ be the circle in ${\mathbb D}$ which has centre $0$ 
and radius $|\omega|^{-n}$, and let $A_{-n}$ be the closed annulus which has boundaries $C_{-n}$ and $C_{-(n+1)}$. \\

Let ${\tilde \Lambda}=\lambda^{-1}(\Lambda)$ (where 
$\Lambda=\Lambda({\mathcal F}_a)=\Lambda_-({\mathcal F}_a)\cup\Lambda_+({\mathcal F}_a)$).
The map $\times \omega$ sends ${\tilde \Lambda}\cap A_{-n}$ bijectively to ${\tilde \Lambda}\cap A_{-(n-1)}$. 
We can extend ${\tilde \Lambda}$ in the obvious way to become a subset 
$\bigcup_0^\infty \omega^n({\tilde \Lambda}) \subset {\mathbb C}$ invariant under $\times \omega$, which 
we also denote by ${\tilde \Lambda}$. We may think of ${\mathbb C}^*={\mathbb C}\setminus\{0\}$ as 
a covering space of the torus obtained from $A_{-n}$ by identifying its boundaries via $\times \omega$.\\

Let $U=U_0$ be any component of ${\mathbb C}^*\setminus {\tilde \Lambda}$. For each $i\in{\mathbb Z}$ 
let  $U_i$ denote the component $\omega^i(U_0)$. Under the map $\times\omega$ we have a 
bi-infinite sequence of components:
$$ \ldots \to U_{-n}\to \ldots \to U_{-1}\to U_0\to U_1\to \ldots \to U_n\to\ldots$$
Intersecting these components with ${\mathbb D}$ and applying the linearisation map $\lambda$
we have a corresponding bi-infinite sequence of maps:
$$f_i: \lambda(U_i\cap {\mathbb D}) \to \lambda(U_{i+1}\cap {\mathbb D})$$
where each $f_i$ is the branch of the correspondence which fixes $\hat{z}$, denoted $f$ a few lines above. 
The reason for the subscript `$i$' is that if we now conjugate our sequence of maps 
$(f_i)_{i\in {\mathbb Z}}$ by the B\"ottcher 
conformal bijection: 
$\varphi=\varphi_a : \widehat \C\setminus \Lambda \to {\mathbb H}$
we obtain a sequence of maps
$$h_i: \varphi(\lambda(U_i\cap {\mathbb D})) \to \varphi(\lambda(U_{i+1}\cap {\mathbb D}))$$
each of which is a restriction of either $\alpha:{\mathbb H} \to {\mathbb H}$ or 
$\beta:{\mathbb H}\to {\mathbb H}$, and it is helpful to have a notation which allows us to distinguish 
these.

\begin{prop}\label{extension}
For each 
$i\in{\mathbb Z}$ the linearising map 
$\lambda:U_i\cap{\mathbb D} \to \widehat \C\setminus\Lambda$ extends to a conformal bijection
$\lambda_i:U_i\to \widehat \C\setminus\Lambda$.
The bijections $\varphi\circ \lambda_i:U_i\to {\mathbb H}$ send the (bi-infinite) sequence of maps `$\times\omega$' 
$$ \ldots \to U_{-n}\to \ldots \to U_{-1}\to U_0\to U_1\to \ldots \to U_n\to\ldots$$
to the corresponding sequence of maps $h_i:{\mathbb H} \to {\mathbb H}$, in a commuting 
ladder.
\end{prop}

{\bf Proof}

First note that each $f_i$ is not just the locally defined branch of the correspondence fixing $\hat{z}$,
but is a well-defined conformal bijection 
$\widehat \C\setminus \Lambda \to \widehat \C\setminus \Lambda$, the conjugate, via $\varphi^{-1}$,
of $h_i:{\mathbb H} \to {\mathbb H}$, that is to say of either $\alpha$ or $\beta$.\\

Now, given any $x_0\in U_i$, choose $n$ sufficiently large that $x_{-n}=\omega^{-n}x_0$ lies in 
$U_{i-n}\cap {\mathbb D}$. Define $\lambda_i(x_0)$ to be:
$$\lambda_i(x_0)=f_{i-1}\circ f_{i-2}\circ\ldots\circ f_{i-n}\circ\lambda(x_{-n}).$$
Proving that $\lambda_i(x_0)$ is well-defined, and that when $x_0\in U_i\cap{\mathbb D}$ the definition 
agrees with that of $\lambda(x_0)$, is a  straightforward exercise. To see that 
$$\lambda_i:U_i \to \widehat \C\setminus \Lambda$$ is injective, observe that if we are given any two 
distinct points in $U_i$ then by applying $\omega^{-n}$ with $n$ sufficiently large we can pull them back
to a pair of  distinct points in $U_{i-n}\cap{\mathbb D}$; this pair then maps forward under the bijection
$f_{i-1}\circ f_{i-2}\circ\ldots\circ f_{i-n}$ to a pair of distinct points in $\widehat \C\setminus \Lambda$.\\

To prove that $\lambda_i$ is surjective, we consider the sequence of conformal bijections
$$f_{i-1}^{-1},\ f_{i-2}^{-1}\circ f_{i-1}^{-1},\ f_{i-3}^{-1}\circ f_{i-2}^{-1}\circ f_{i-1}^{-1},\ldots$$  
from $\widehat \C\setminus \Lambda$ to itself. \\

This sequence forms a normal family (since $\widehat \C\setminus \Lambda$ is a hyperbolic surface), 
so some subsequence converges locally uniformly to a holomorphic map 
$\widehat \C\setminus \Lambda \to \widehat \C$. But on any compact subset of
$\lambda(U_i\cap {\mathbb D})$ the whole sequence converges uniformly to the map which
sends the compact set to the fixed point $\hat{z}$. By uniqueness of analytic continuation it follows that the 
only holomorphic map to which any subsequence of our self-maps of $\widehat \C\setminus \Lambda$
can converge is the constant map $$\widehat \C\setminus \Lambda \to \hat{z}\in \widehat \C.$$
So for every $z \in \widehat \C\setminus \Lambda$ the images
 $$f_{i-1}^{-1}(z),\ f_{i-2}^{-1}\circ f_{i-1}^{-1}(z),\ f_{i-3}^{-1}\circ f_{i-2}^{-1}\circ f_{i-1}^{-1}(z),\ldots$$
converge to $\hat{z}$. Thus for some $n$ the image 
$f_{i-n}^{-1}\circ f_{i-n+1}^{-1}\circ\ldots\circ f_{i-1}^{-1}(z)$ of $z$ lies in
$\lambda(U_j\cap {\mathbb D})$ for some $j$. But, since the images of a line segment joining $z$
to a point of $U_i\cap {\mathbb D}$ under the sequence of maps must also converge uniformly to the single 
point $\hat{z}$, it is easily seen that $j=i-n$. It follows from our definition of the extension $\lambda_i$ that 
$z\in \lambda_i(U_i)$. 
That the homeomorphisms $\lambda_i$ form a commuting ladder follows from their definition, since they 
form such a ladder when restricted to the $U_i\cap {\mathbb D}$.
\qed
\\

For a repelling cycle $\{z_0,\ldots z_{m-1}\}$, $m>1$, in place of a fixed point of ${\mathcal F}_a$,
an analogue of the analysis above goes through in the obvious way.\\

Our next observation is there are restrictions on the itineraries that can occur for Fatou components of the linearised map, either 
at a repelling fixed point or a repelling cycle.

\begin{lemma}\label{no_long_sequences}
For each repelling cycle $\{z_0,\ldots z_{m-1}\}$, $m\ge 1$, of ${\mathcal F}_a$, there is a bound on the length of a sequence of 
consecutive occurrences of 
$\alpha$, or consecutive occurrences of $\beta$, that can occur in the itinerary $S$ of a Fatou component $U$ of the linearisation in
a neighbourhood of the cycle.
\end{lemma}

{\bf Proof}

First suppose that $U$ is a Fatou component at the point $z_0$ of the cycle, and the itinerary $S$ of $U$ contains 
$n$ consecutive $\alpha$'s. Replacing $U$ by its appropriate forward or backward
image, we may suppose these $n$ consecutive $\alpha$'s lie immediately to the left of the marker point. Since $\alpha$ 
is the map $z\to z+1$ on ${\mathbb H}$ and its boundary, it follows from Proposition \ref{extension} that we can 
choose a point $x_n\in U$ such that $\lambda(x_n)$ is arbitrarily close to $z_0$, its image $\varphi\lambda(x_n)$ is arbitrarily close 
to the boundary of ${\mathbb H}$, and $Re(\varphi\lambda(x_n))<-n$. \\
 
Now, if for each positive integer $n$ the itinerary $S$ contains $n$ consecutive $\alpha$'s, we may construct a sequence of points $x_n$
in the appropriate components of the linearised map with the properties that

(i) the sequence $(\lambda(x_n))_{n\ge 1}\subset \Omega({\mathcal F}_a)$ converges to $z_0\in \widehat \C$;

(ii) its image under $\varphi$, the sequence $(\varphi\lambda(x_n))_{n\ge 1}\subset {\mathbb H}$, converges to 
$-\infty\in \partial{\mathbb H}$.

But (i) and (ii) are contradictory, since the inverse of the B\"ottcher map extends continuously to $-\infty$, sending $-\infty$
to $P$, the  parabolic fixed point, yet $z_0\ne P$.\\  

We obtain a similar contradiction when we take $\beta$ in place of $\alpha$, and $0$ in 
place of $-\infty$. 
 \qed\\

{\bf Notation.} 
Recall that $\gamma(S)$ is our notation for the geodesic in ${\mathbb H}$ which has itinerary $S$. It is convenient for the 
remaining proofs in the current section to introduce a notation for the image of $\gamma(S)$ in $\Omega({\mathcal F}_a)$
under the inverse $\varphi_a^{-1}$ of the B\"ottcher map. We define 
$$g(S):=\varphi_a^{-1}(\gamma(S))\subset \Omega({\mathcal F}_a)=\widehat \C\setminus \Lambda({\mathcal F}_a).$$  

In a linearising neighbourhood of a repelling fixed point $\widehat \C$ we can further pull back $g(S)$
to $\lambda^{-1}(g(S)) \subset {\mathbb D}$.\\

If $\gamma(S)$ is a geodesic in ${\mathbb H}$ with initial end point  in ${\mathbb R}^-\subset\partial{\mathbb H}$, we parametrise $\gamma(S)$ in a neighbourhood of
this end point by the potential $\chi$ (which we recall is approximately `height') and thereby also parametrise the corresponding part of $g(S)$.
We shall write $g_S$ for the parametrising function.

\begin{lemma}\label{convergent_geodesics} If $\gamma(S_n)$, $n>0$, is any sequence of geodesics which have initial and final end points 
converging to the initial and final end points of $\gamma(S)$, then for each sufficiently small pair of positive real numbers $t_*<t^*$ the
geodesics $g(S_n)$ converge uniformly to $g(S)$ on the interval $[t_*,t^*]$. 
\end{lemma}

{\bf Proof}

Working in the disc model of hyperbolic space, the inverse of the 
B\"ottcher map is uniformly continuous in every closed annulus centred at the centre of the disc. Since the geodesics $\gamma(S_n)$
converge uniformly to the geodesic $\gamma(S)$ on the compact set $[t_*,t^*]$, the result follows. \qed\\

Next, in place of the `fundamental domain' $I_t(g_s)$ on the `dynamical ray' $g_s$, defined in \cite{BLyu}, we define
a `basic interval' on the geodesic $g(S)$:

\begin{defn}
The basic interval is the  subset $I_t(g(S)):= g_S[t/\lambda_+,t]$ of $g(S)$. 
\end{defn}

Here $\lambda_+$ is an upper bound on the multiplier for $\alpha$ and $\beta$ on the potential of a point 
$z$ close to the 
negative real axis in ${\mathbb H}$, as computed in Lemma \ref{alphabetapotential}(i) (Section \ref{Green_fn}).
As $t$ tends to $0$  the hyperbolic length of $I_t(g_S)$ tends to $\log(\lambda_+)$. Recall that $x^-(S)$ denotes the
landing point of the geodesic $\gamma(S)\subset {\mathbb H}$ on ${\mathbb R}^-$. 

\begin{prop}\label{basic_length_to_zero} For each $K>1$, the Euclidean length of every basic interval 
$I_t(g(S))$ with $x^-(S)\in [-K,-1/K]$ tends to zero, independently of $S$, as $t$ tends to $0$.
\end{prop}

{\bf Proof} 

Approaching the boundary of any bounded simply-connected domain in the plane, the density of a hyperbolic metric on the 
domain
tends to zero (for a proof, see for example Lemma 2.3 of \cite{BLyu}). For small $t>0$ 
the basic intervals $I_t(g(S))$ with $x^-(S)\in [-K,-1/K]$ have uniformly bounded hyperbolic lengths and are contained in 
the neighbourhood of $\partial \Omega({\mathcal F}_a)$ bounded by the equipotential $\{z\in \Omega: G(z)=t\}$ (where $G$ 
is the potential function
we chose in Definition \ref{Green} in Section \ref{Green_fn}). Hence the Euclidean length of $I_t(g(S))$ tends to zero
with $t$. Uniformity with respect to $S$ follows from the compactness of each equipotential. \qed\\

We now have all the ingredients to prove that there is at least one geodesic landing at each repelling periodic point, inspired by the methods in \cite{BLyu}.\\

{\bf Proof of Theorem \ref{repelling_fp}}\\

The proof is divided into four steps. For clarity of exposition, in the first three steps we consider the case of a repelling fixed point 
${\hat z}$ and afterwards in Step 4 we list the modifications to Steps 1, 2 and 3 needed to prove the result for a repelling cycle of 
period $m>1$.\\

{\it Step 1: there exists a landing geodesic}\\

We follow the strategy of the proof of Theorem 2.5 of \cite{BLyu}, constructing a landing geodesic as a limit of
longer and longer segments of a convergent sequence of geodesics,
the major difference from \cite{BLyu} being that we add an overlapping `basic interval' of geodesic at each stage, 
rather than adding a `fundamental domain' to the end of the segment of geodesic already constructed. \\

Let $U$ be a component of the Fatou set in a linearising neighbourhood  around $\hat{z}$
and $S$ be its itinerary. If $S$ is periodic, then so is $g(S)$, and in this case $g(S)$ lands, by  Proposition \ref{periodic_rays_land}. 
Otherwise, write $g_i$ for the $-i$th shift of $g_0=g(S)$ (so $g_i=f_a^{-i}(g_0)$), and write $\gamma_i$ for the 
geodesic $\varphi_a(g_i) \subset {\mathbb H}$. Near to the boundary of ${\mathbb H}$ the geodesic 
$\gamma_i$ is parametrised by the potential function $\chi$, and near to $\Lambda$ the geodesic $g_i$ is parametrised
by $G=\chi\circ\varphi_a$. We note that by Lemma \ref{no_long_sequences} there is a bound on the length 
of consecutive appearances of the same letter in $S$, and it follows that the initial points of the geodesics $\gamma_i$ 
fall within some closed interval $[-K,-1/K]$ contained within $(-\infty,0)\subset\partial{\mathbb H}$. So, by
Lemma \ref{alphabetapotential}(ii) (Section \ref{Green_fn}), there exists a `minimum multiplier' $\lambda_-$ for the effect of 
$\alpha$ and $\beta$ on potential.\\

Let $V'$ be a linearising neighbourhood around $\hat{z}$, and $V\subset V'$ be
its pre-image $f_a^{-1}(V')$. Let $\varepsilon$ be the (Euclidean) distance between $\partial V$ and $\partial V'$.
Let $I_t(g_i)$ denote $\varphi_a^{-1}(I_t(\gamma_i))$.
We note that the Euclidean length of $I_t(g_i)$  tends to zero (uniformly in $i$) as $t$ tends to $0$ by Proposition \ref{basic_length_to_zero}.
Thus there exists $t_\varepsilon>0$ such that the Euclidean length of 
$I_t(g_{i})$ is less than $\varepsilon$ for all $i\ge 0$ and $t<t _\varepsilon$.  \\

We claim there exists some $i$ and $t_0<t_\varepsilon$ with $z_0:=g_i(t_0)\in V'$. To see this, observe that 
$\widehat \C\setminus\Lambda$ is foliated by images under $\varphi_a\circ\lambda$ of circles in ${\mathbb D}$, and at least one
of these images of circles must meet $\gamma_i$. Now pulling back by $f_a^{-1}$
pulls points in the sphere back by $\omega^{-1}$ (so we can pull back into $V$)
and reduces potential in ${\mathbb H}$ by a factor of at least $\lambda_-$ (so we can ensure $t_0<t_\varepsilon$).
Renumber this $g_i$ as $g_0$  and renumber the $g_n$ accordingly, so that $g_n=f_a^{-n}(g_0)$.\\

\begin{figure}
\begin{center}
\scalebox{.45}{\includegraphics{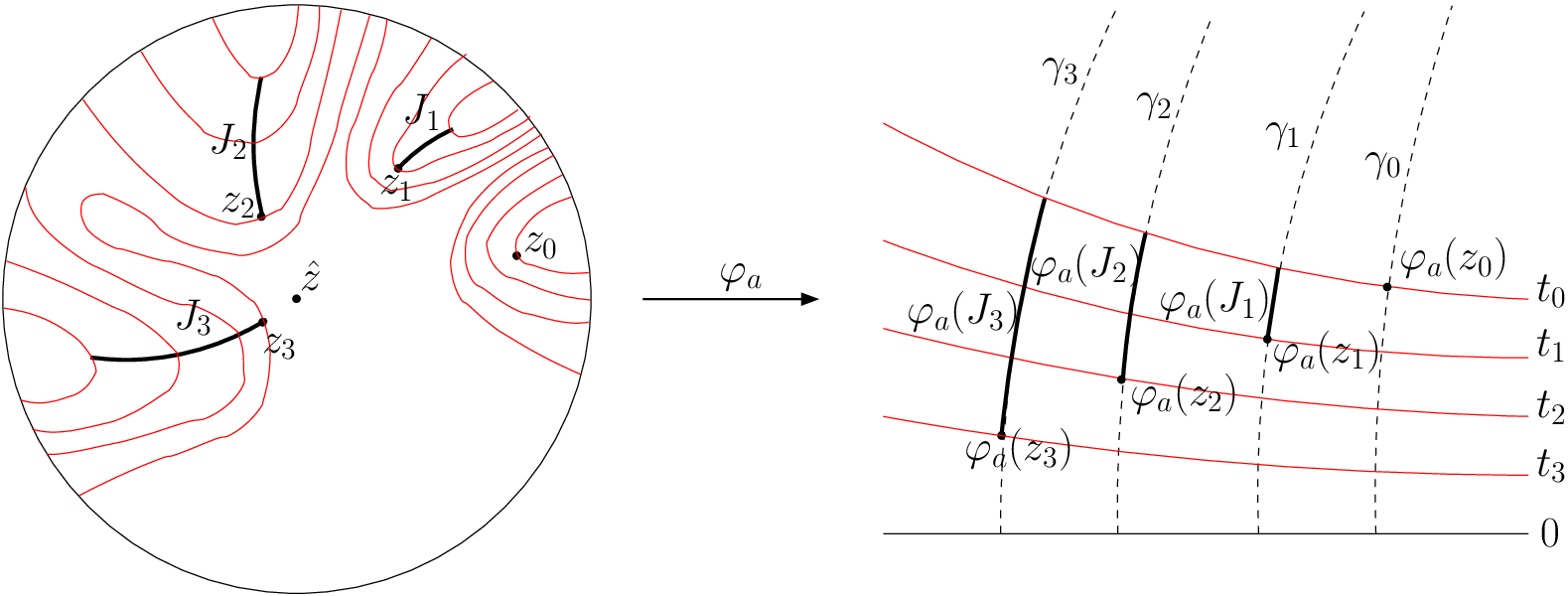}} 
\caption{\small Step 1 of the proof of Theorem 2. On the left, the neighbourhood
$V'$ of $\hat{z}$, and on the right the upper half-plane ${\mathbb H}$. The B\"ottcher map $\varphi_a$ sends $V'\cap 
(\widehat \C\setminus\Lambda)$ injectively into ${\mathbb H}$.}\label{thm2fig}
\end{center}
\end{figure}

Let 
$z_n=f_a^{-n}(z_0)\in g_n$, let $t_n=\chi(\varphi_a(z_n))$, let $J_1=g_1[t_1,t_0]$ and, inductively for $n>1$,
let $J_n=f_a^{-1}(J_{n-1})\cup I_{t_0}(g_n)$.
Figure \ref{thm2fig} illustrates the points $z_n$ and geodesic segments $J_n$ together 
with their $\varphi_a$-images in  ${\mathbb H}$.
For each $n\ge 1$, $\varphi_a(z_n)=h_n^{-1}(\varphi_a(z_{n-1}))$, where $h_n$ is $\alpha$ or $\beta$ 
(as defined in Proposition \ref{extension}). 
As $h_n^{-1}(\varphi_a(J_{n-1}))$ overlaps 
$I_{t_0}(\gamma_n)=\gamma_n[t_0/\lambda_+,t_0]$ (by Lemma \ref{alphabetapotential}(ii), Section \ref{Green_fn}), 
we deduce that $\varphi_a(J_n)\supseteq\gamma_n[t_n,t_0]$.\\

The point $z_1=f_a^{-1}(z_0)\in V$, so $J_1\subset V'$ (since the Euclidean length of $I_{t_0}(g(S_1))<\varepsilon$). Inductively,
by similar reasoning, each of the intervals $J_n\subset g(S_n)$ lies in $V'$. Moreover the sequence 
$(t_n)$ of lower ends of these intervals has limit zero (since, inductively, $t_0/t_n>(\lambda_-)^n$).\\

Within our sequence $(g_n)$ of geodesics in $\Omega({\mathcal F}_a)$ 
we may choose a subsequence $(g_{n_j})$ such that
the left hand ends of the corresponding $\gamma_{n_j}$ in ${\mathbb H}$ form a convergent sequence, and their
right hand ends also converge. The fact that the itinerary $S$ has a bound on the length of consecutive repetitions of 
the same symbol 
ensures that the left hand ends of the geodesics $\gamma_n$ are bounded away 
from $0$ and $\infty$ and hence that the $\gamma_{n_j}$ converge to a genuine geodesic, not to the `degenerate geodesic' consisting of
the single point $0$ or the single
point $\infty$. Thus the sequence of geodesic segments $(J_{n_j})$ converges to a segment $g(S')(0,t_0)$ 
of some limit geodesic $g(S')$. For any $t_*$ such that $0<t_*<t_0$ this convergence is uniform on the interval $[t_*,t_0]\subset (0,t_0)$ 
by Lemma \ref{convergent_geodesics}.
Using the fact that close to the 
fixed point $\hat{z}$ the map $f_a$ is approximated by $z \to {\hat z}+\omega(z-{\hat z})$, it can now be verified easily that 
$g(S')(0,t_{n_j}]\subset B({\hat z}, A/|\omega|^{n_j})$ for some constant $A$ and $n_j$ sufficiently large, and thus that
$g(S')$ restricted to $(0,t_0]$ can be extended continuously to $t=0$ by setting $g(S')(0)=\hat{z}$.
Then $g(S')$ is a 
geodesic landing at the fixed point.\\

{\it Step 2: every landing geodesic is periodic}\\

Given any sequence $(\gamma_n)$ of geodesics in ${\mathbb H}$ with left hand ends $x^-(\gamma_n)$ converging to a point 
$x\in {\mathbb R}^{<0}$, and with the corresponding geodesics $g_n=\varphi_a^{-1}(\gamma_n)$ in $\Omega({\mathcal F}_a)$ 
all landing at $\hat{z}$, the 
method of Step 1 can be applied to construct a geodesic $\gamma$ in ${\mathbb H}$ with $x^-(\gamma)=x$ and 
with $\varphi_a^{-1}(\gamma)$ also landing at $\hat{z}$. Thus the set $A$ of left-hand end-points of images under $\varphi$ of geodesics in
$\Omega({\mathcal F}_a)$ which land at $\hat{z}$ is a {\it closed} subset of ${\mathbb R}^{<0}$. This set $A$, regarded as a subset of 
the circle obtained by identifying the end points $-\infty$ and $0$ of ${\mathbb R}^{<0}$, is invariant under the doubling map defined
by $\alpha$ on $(-\infty,-1)$ and $\beta$ on $(-1,0)$, and its cyclic order is preserved by this doubling map. 
Thus $A$ is a closed invariant Sturmian subset of ${\mathbb R}/{\mathbb Z}$. Moreover the doubling map is injective on $A$, since $\F_a$ is
linearisable in a neighbourhood of $\hat z$. By Corollary \ref{invariant_sets_finite} (Section \ref{St}), it follows that
$A$ is the (unique) finite Sturmian orbit of rotation number $p/q$ for some rational $p/q$.\\

{\it Step 3: counting cycles of landing geodesics}\\

By Step 2 the itinerary of a landing geodesic is necessarily  of
the form $W^\infty$ where $W=W_{p/q}$ is the Sturmian word of rotation number $p/q$. There is only one such orbit.  \\

{\it Step 4: modifications to Steps 1 to 3 needed to prove the Theorem for a repelling cycle of period $m>1$}\\

Step 1 is unchanged except that $\hat{z}$ is replaced by one of the points of the repelling cycle and $f_a$ is replaced by the 
first return map $f_a^m$.\\ 

In Step 2 in the case $m>1$, the closure of the set of left-hand end-points of $\varphi_a$-images of geodesics landing
at points of the cycle becomes a union of $m$ disjoint closed sets, each of which is invariant under the appropriate word $W$ of length 
$m$ in the symbols $\alpha$ and $\beta$ (the first return map) and has its cyclic order 
preserved by $W$, and each of which has the same well-defined rotation number $\nu$ under the first return map. 
To conclude the proof of this step we
must exclude the possibility that these sets are infinite. Rather than generalising Corollary \ref{invariant_sets_finite} to
this non-Sturmian situation, the easiest way to proceed is to apply the even more general 
Lemma 2.6 of \cite{BLyu}, which states that if $f$ is a locally expanding map of a compact metric 
space $X$ to itself, and $A$ is a closed invariant subset of $X$ restricted to which $f$ is invertible, then $A$ is finite.
(Of course this lemma also provides an alternative proof of Step 2 for a fixed point.) \\

For Step 3 in the case $m>1$, we observe that the combinatorics of itineraries of
geodesics landing on repelling cycles for correspondences in our family are identical to the combinatorics 
of rays landing on repelling cycles for quadratic polynomials 
(see Milnor \cite{M1}, or Schleicher \cite{S}, for the latter: each step in these analyses goes through in the same 
way for our correspondences 
${\mathcal F}_a$). The fact that there are at most two orbits of rays which 
land on any particular
repelling cycle for a quadratic polynomial follows from the `tuning' theory of Douady and Hubbard \cite{DH}. 
The two orbit case corresponds to {\it primitive} components of the 
Mandelbrot set $\mathcal M$: at the root point $c$ of such a component the quadratic map $Q_c:z\to z^2+c$ has a parabolic cycle 
which has exactly two orbits of 
landing rays, and this landing pattern persists when the parabolic orbit is deformed into a repelling orbit. The single orbit case 
corresponds to landing patterns born at the root points of {\it satellite} components of $\mathcal M$:  the landing rays in such patterns have 
binary sequences which are renormalisable.
\qed \\

Finally we note in the corollary below that for a repelling fixed point $\hat{z}$ not only there is a unique cycle of periodic geodesics 
which land at $\hat{z}$, but there is just one 
cycle of Fatou components of the linearisation there. Similarly it can be proved that in the case of a repelling cycle of
period $m>1$ the linearised map has either exactly one or exactly two cycles of Fatou components.

\begin{cor} The Fatou components of the linearised map at a repelling fixed point ${\hat z}$ form a single cycle.
\end{cor}

{\bf Proof} 

The periodic geodesic landing at $\hat{z}$ has itinerary $W^\infty$ where $W=W_{p/q}$ for some $p/q$. If $U$ and $V$ are Fatou 
components of the linearisation then both must also have itinerary $W^\infty$. But now
$U$ and $V$ must be the same component, for we can join any point of $\lambda(U$) to a point of $\lambda(V)$ by a path in 
$\Omega$, and now applying $W^{-n}$ for sufficiently large $n$ shrinks this to a path contained within an arbitrarily small neighbourhood
of ${\hat z}$. \qed\\

\section{The proof of Theorem \ref{Yoccoz}}\label{Yoccoz_proof}
We recall \cite{DH} that a quadratic map $z\to z^2+c$ with $c\in \mathcal M$ has two fixed points: the 
beta-fixed-point is the landing point of the external ray of argument zero, so it has combinatorial rotation number $0$.
The other fixed point, known as the 
alpha-fixed-point, is a repeller precisely for those $c\in  \mathcal M$ which lie outside the closure of the main cardioid. 
Generically a correspondence in the family $\F_a$ has $4$ fixed points,  
the parabolic fixed point at $\Lambda_-\cap\Lambda_+$  ($Z=1$, or equivalently, in the $z$-coordinate, $z=0$), which
is always a double fixed point, and two others. When $a \in {\mathcal K}$, the Klein combination locus, these two others are
one each in $\Lambda_-$ and $\Lambda_+$. We are concerned here with the fixed point in $\Lambda_-$, call it $p_a$, in the 
Douady-Hubbard terminology the ``alpha-fixed-point'' of 
the $2$-to-$1$ branch $f_a$ of ${\mathcal F}_a$ defined earlier. Our Yoccoz inequality establishes bounds on the derivative $f_a'(p_a)=\zeta_a$, 
when $a \in {\mathcal M}_\Gamma$ and $p_a$ is repelling.\\

Let $U_0$ be a (periodic) Fatou component of the linearised map at the repelling fixed point $p_a \in \Lambda_{a,-}$ 
(as in the preceding section), and let 
${\mathbb T}$ denote the torus ${\mathbb C}^*/(\times \zeta)$, so ${\mathbb T}$ is the quotient of ${\mathbb C}$
by the lattice generated by $z \to z+2\pi i$ and $z\to z+\tau$, where $\tau$ is the principal value of the complex logarithm of $\zeta$.\\

The image of $U_0$ in ${\mathbb T}$ 
is an annulus ${\mathcal A}$, homotopic to a closed $(-p,q)$-curve wrapping around ${\mathbb T}$.
Since this annulus is embedded in ${\mathbb T}$, we have, for a suitable choice of $\tau$ mod $2\pi i$,
$$mod({\mathcal A})\le \frac{2\pi Re(\tau)}{|2 \pi i p - \tau q|^2}$$
(see \cite{H}, Proposition 3.2, for a justification).
However it follows from Proposition \ref{extension} (Section \ref{repellor}) that ${\mathcal A}=U_0/(\times \zeta^q)$ maps 
bijectively to the annulus ${\mathbb H}/T_{p/q}$,
where $T_{p/q}$ is the Sturmian word of length $q$ in the letters $\alpha$ and $\beta$ which corresponds to
rotation number $p/q$. But $T_{p/q}$ acts on ${\mathbb H}$ by
$$z \to \mu(T_{p/q})z$$
and so the annulus ${\mathbb H}/T_{p/q}$ is that obtained from the region $1\le |z|\le \mu(T_{p/q})$ in ${\mathbb H}$ by identifying
the bounding semicircles. Mapping $z$ to $\log{z}$ sends this region to a rectangle of side lengths $\pi$ and 
$\log{(\mu(T_{p/q})}$, so ${\mathbb H}/T_{p/q}$ has modulus $\pi/\log{\mu(T_{p/q})}$.
By applying Corollary \ref{estimates} (Section \ref{St}) we deduce that in the case $p/q\le 1/2$:
$$mod({\mathcal A})>\frac{\pi}{2p\log(\lceil q/p \rceil+1)}.$$
From the two inequalities above, we have (still in the case $p/q\le 1/2$):
$$\frac{Re(\tau)}{|\tau-2\pi i p/q|^2}\ge \frac{q^2}{4p\log(\lceil q/p \rceil+1)},$$
which is equivalent to the statement for $p/q\le 1/2$ in the theorem. (Given any $r>0$, the set of points $z=x+iy$ which 
satisfy the inequality $x\ge |z|^2/2r$ form a disc of radius $r$ tangent to the imaginary axis at the origin.) 
The statement for $p/q\ge 1/2$ follows from Corollary \ref{estimates} in the same way. \qed\\

From Theorem \ref{Yoccoz} we obtain the following practical criterion:

\begin{cor}\label{Yoccoz_cor}
Let $a\in {\mathcal C}_\Gamma$. If the derivative $\zeta$ of $f_a$ at a repelling fixed point 
has its argument in the interval $(0,\pi]$, then the principal value of the complex logarithm of $\zeta$ 
$$\tau=\log(|\zeta|)+iArg(\zeta)$$  
lies in the part of 
${\mathbb H}$ defined by 
$$\Re(\tau)<4.3\nu^2\log(\nu^{-1}+1) \  {where} \ \nu=\Im(\tau)/2\pi.$$
\end{cor}

Here the multiplying factor of $4.3$ is chosen to ensure that the curve lies outside the 
union of the discs permitted by our Yoccoz inequality, not only outside their horizontal diameters. 
Note that as the graph is concave it suffices to find a multiplying factor such that the tangent to
the curve where it crosses the horizontal $y=\pi$ does not meet
$D_{1/2}$, the disc corresponding to $p/q=1/2$.
However we remark that we can sharpen our estimates of disc radii 
for rotation numbers of the form $\nu=1/q$, by computing the multiplier of $\alpha^{q-1}\beta$ exactly.
The disc $D_{1/2}$ corresponds to the Sturmian word $\alpha\beta$, which has multiplier 
$((3+\sqrt{5})/2)^2$ (as we saw in the proof of Proposition \ref{multipliers}). 
Repeating the calculation in the proof of Theorem \ref{Yoccoz} for the case $p/q=1/2$, but now replacing the estimate 
$(\lceil 2/1 \rceil+1)=3$
by the sharper value $(3+\sqrt{5})/2$ we obtain the value $\log((3+\sqrt{5})/2)$ for the diameter of $D_{1/2}$, the largest
disc, and hence the following 
absolute bound on the modulus of the derivative:

\begin{cor}\label{abs_bound}
If $a\in {\mathcal C}_\Gamma$ then the derivative $\zeta$ of $f_a$ at its $\alpha$-fixed-point
satisfies the inequality $$|\zeta|\le \frac{3+\sqrt{5}}{2}.$$
\end{cor}

In fact this bound is sharp (see Remark \ref{bound} in Section \ref{param_lune} below).\\

In Figure \ref{Yoccoz_circles} in the Introduction, we plot some of the discs in the $\log{\zeta}$-plane permitted by the 
Yoccoz inequality, on the left for 
matings between quadratic polynomials and $PSL(2,{\mathbb Z})$, and on the right for quadratic polynomials (here
the disc for $\nu=p/q$ had radius $(\log{2})/q$). In the left-hand picture the 
discs lie entirely to the left of the curve (also illustrated):
 $$\{\tau=4.3\nu^2\log(\nu^{-1}+1)+2\pi \nu i:\ 0\le \nu\le 1/2\}.$$ 

\begin{remark} $\Re(\tau)$ has faster convergence to $0$ as $\Im(\tau)\to 0$ than is the case for the classical Yoccoz 
inequality for quadratic polynomials, where the corresponding region in the upper half-plane is bounded by a 
straight line (see Figure \ref{Yoccoz_circles}).
The underlying reason for the linear bound in the classical case (the right-hand picture in Figure \ref{Yoccoz_circles}) is 
that for the map $z \to 2z$ on $\H/\Z$ the multiplier of {\it every} period $q$ cycle is $2^q$. However in our case 
we are dealing with cycles made up from the parabolic maps $\alpha$ and $\beta$ on $\H$, and the multiplier of a $q$-cycle varies with the word 
$W$ which defines the cycle. The curved bound in the left-hand picture arises from the bounds we computed for the multipliers of the cycles $\alpha^{q-1}\beta$
corresponding to rotation number $1/q$ (Proposition \ref{multipliers}(i)): these multipliers grow {\it quadratically} in $q$, not {\it exponentially}.
\end{remark}

\begin{remark}\label{higher_period} An obvious question is the nature of the Yoccoz inequalities for repelling periodic orbits of period greater than $1$. As an example we compute these for period $2$ orbits in the case that the rotation number
of the first return map is of the form $1/q$. Tuning theory (renormalisation) tells us that the itinerary of such a geodesic landing on a 
period two cycle is (a cyclic permutation of) $(W_{1/q})^\infty$ where $W_{1/q}=BA^{q-1}$ with $A=\alpha\beta$ and $B=\beta\alpha$. 
We easily compute that
$$tr(W_{1/q})=Fib_{2q+1}+Fib_{2q-3}$$
where $Fib_n$ denotes the $n$th Fibonacci number ($F_{1}=1,\ F_{2}=1,\ F_3=2, \ldots$).
Since $Fib_n=((1+\sqrt{5})/2)^n - (1-\sqrt{5})/2)^n)/\sqrt{5}$ we can compute an exact formula for $tr(W_{1/q})$, and as 
$tr(W_{1/q})=\mu_{1/q}+\mu^{-1}_{1/q}$ where $\mu_{1/q}$ is the multiplier of $W_{1/q}$ at its fixed point, we may deduce that, 
for large $q$,
$$\mu_{1/q}\sim \frac{1}{5}
\left(\frac{1+\sqrt{5}}{2} \right)^{2(2q+1)}.$$
Thus for a repelling period $2$ orbit of a correspondence in the 
family ${\mathcal F}_a$, having first return map of rotation number $1/q$, 
the analogue of Figure \ref{Yoccoz_circles} gives a linear bound, as in the quadratic
polynomials case.
\end{remark}

\section{Proof of Theorem \ref{param-lune} and Corollary \ref{near_7}}\label{param_lune}
We first compute the derivative of $f_a$ at its alpha-fixed-point, in preparation for applying our 
Pommerenke-Levin-Yoccoz inequality to prove 
Theorem \ref{param-lune}. \\

We use the $z$ coordinate system, where the correspondence ${\mathcal F}_a$ has equation:
$$\left(\frac{az+1}{z+1}\right)^2+\left(\frac{az+1}{z+1}\right)\left(\frac{aw-1}{aw+1}\right)+\left(\frac{aw-1}{w-1}\right)^2=3. \leqno(*)$$
To find the alpha-fixed-point $z_0$ we set $w=z=z_0$ in $(*)$ and get
$$(3a^2-3)z_0^4+(a^2-8a+7)z_0^2=0$$
which, ignoring the parabolic fixed point $z_0=0$ (the beta-fixed-point), gives us
$$z_0=-\sqrt{\frac{7-a}{3(a+1)}}$$
where the minus sign before the square root indicate that we are taking the branch which has a negative value 
when $a$ is real and between $+1$ and $7$.\\

To find the derivative at $z_0$ we differentiate $(*)$ with respect to $z$, and set $z=w=z_0$. We deduce 
that at $z=z_0$ the value of $dw/dz$ is
$$\zeta= \left(\frac{z_0-1}{z_0+1}\right)^2 \left(\frac{3(az_0^2-1)+(1-a)z_0}{3(az_0^2-1)+(a-1)z_0}\right).$$
After we substitute $z_0^2=(7-a)/[3(a+1)]$, 
$$\zeta=\frac{(a^2-2a-11)+(a+1)(7-a)z_0}{(a^2-2a-11)-(a+1)(7-a)z_0}.$$

\begin{remark}
Notice that the expression for $\zeta$ can be written in the form $$\zeta=\frac{1+E}{1-E}$$
and that $|\zeta|=1$ if and only if $E$ is pure imaginary; also that $\zeta$ is real if and only if 
$E$ is real.
Values where $|\zeta|=1$ coincide with the boundary of the
main component of the interior of $\mathcal{M}_{\Gamma}$. For example the value $\zeta=-1$,  where the boundary 
of the main component 
cuts the real axis, is given by the positive solution to $a^2-2a-11=0$, that is $a=1+2\sqrt{3}=4.464$. 
\end{remark}

Denote the open disc $\{a:|a-a_0|<r\}$ by $D(a_0,r)$. 
To prove Theorem \ref{param-lune} it will suffice to prove the following three statements:\\

(i) There exists $\delta>0$ such that $\mathcal{M}_{\Gamma}$ does not meet $D_1=D(1,\delta)$; \\

(ii) For every $\alpha$ such that $\pi/3<\alpha\le\pi/2$, there is a disc neighbourhood $D_7=D(7,\epsilon)$
with centre $a=7$ and radius $\epsilon>0$ such that $\mathcal{M}_{\Gamma}\cap D_7\subset L_\alpha\cap D_7$;\\

(iii) $\partial{D(4,3)}\setminus (D_1\cup D_7)$ has a neighbourhood which does not meet $\mathcal{M}_{\Gamma}$.\\

{\bf Proof of (i).} This is immediate from our formula for $\zeta$, which gives $\zeta \to \infty$ as $a \to 1$, so our Yoccoz inequality is violated for
$a$ in a disc $D(1,\delta)$.\\

For the proof of (ii) and (iii),
we change the parameter. We first note that since the lower half of the lune boundary is complex conjugate to the 
upper half, it will suffice to prove (ii) and (iii) for points on the upper half. Now let
$$b=\frac{a-7}{a-1}.$$
The upper half of the lune boundary becomes the straight line $$b=te^{i(\pi-\alpha)}, \quad t \in{\mathbb R}^{\ge 0},$$
the point $z_0$ becomes
$$z_0=-\sqrt{\frac{b}{b-4}},$$
and $\zeta=\zeta(a)$ becomes
$$\zeta(b)=\frac{2+2b-b^2+b(b-4)z_0}{2+2b-b^2-b(b-4)z_0}$$
which we may write as
$$\frac{1+E}{1-E}$$
where
$$E=\frac{b(4-b)}{2+2b-b^2}\sqrt{\frac{b}{b-4}}.$$

{\bf Proof of (ii).} Substituting $b=te^{i(\pi-\alpha)}$ in our expression for $E$ gives
$$E=\frac{-te^{-i\alpha}(4+te^{-i\alpha})}{2-2te^{-i\alpha}+t^2te^{-2i\alpha}}\sqrt{\frac{te^{-i\alpha}}{te^{-i\alpha}+4}},$$
the leading term of which is
$$-t^{3/2}e^{-3i\alpha/2},$$
and so  we deduce that as $t\to 0$
$$\log{\zeta}=\log{\left(\frac{1+E}{1-E}\right)}\sim 2E\sim -2t^{3/2}e^{-3i\alpha/2}=2t^{3/2}e^{i(\pi-3\alpha/2)}.$$
Thus as $t\to 0$ the complex number $\log{\zeta}$ approaches $0$ tangentially to a straight line 
of argument $\pi-3\alpha/2$. When $\pi/3\le \alpha \le \pi/2$ this line is in the positive quadrant:
when $\alpha=\pi/3$ it is the imaginary axis and when $\alpha=\pi/2$
it is the line of argument $\pi/4$. By Corollary \ref{Yoccoz_cor} it follows that for every value of $\alpha$ in 
the interval $\pi/3<\alpha\le\pi/2$ the set $\mathcal{M}_{\Gamma}\cap D(7,\epsilon)$ lies inside $L_\alpha\cap D(7,\epsilon)$
for a sufficiently small value of $\epsilon>0$.\\

{\bf Proof of (iii).} It will suffice to show that the entire parameter arc $b=it$, $0<t<\infty$, lies outside the region 
permitted by the Yoccoz inequality, since by continuity this will imply that the intersection of the arc with the 
complement of discs centred at its two ends has a neighbourhood which misses $\mathcal{M}_{\Gamma}$. We consider 
$0<t \le 1$ and $1<t<\infty$ separately. \\

For $b=it$ we have $$E=\frac{it(4-it)}{2+2it+t^2}\sqrt{\frac{-it}{4-it}}=\frac{t^2+4it}{(2+t^2)+2it}\sqrt{\frac{-it}{4-it}}.$$
Now 
$$\left|\frac{t(4i+t)}{(2+t^2+2it}\right|^2=t^2\left(\frac{16+t^2}{4+8t^2+t^4}\right)$$
which is strictly increasing for $t\in [0,1]$ so has value $<\sqrt{17/13}<6/5$ there, and for $0\le t \le 1$ we also have
$$\left|\sqrt{\frac{-it}{4-it}}\right|<1/2$$
and so for this range of $t$ we have
$$|E|<\frac{3}{5}.$$
We next estimate the arguments of the two factors of $E$.
$$\frac{t^2+4it}{(2+t^2)+2it}=\frac{t(4i+t)((2+t^2)-2it)}{(2+t^2)^2+4t^2}=\frac{t(t^3+8t^2+2t+(8+2t^2)i)}{(2+t^2)^2+4t^2}$$
which certainly has argument in $[0,\pi/2]$ for all $t\in [0,1]$.\\

Also 
$$\arg{\left(\frac{-i}{4-i}\right)}=\arg{(1-4i)}=-\arctan{4}$$
and as $\arctan{4}>5\pi/12$ we deduce that 
$$\arg(E)<\frac{\pi}{2}-\frac{5\pi}{24}=\frac{7\pi}{24}.$$
Our formula for the derivative $\zeta$ at $z_0$ is
$$\log(\zeta)=\log\left(\frac{1+E}{1-E}\right)=2E(1+E^2/3+E^4/5+E^6/7+\ldots)$$
so for $|E|<3/5$ we have
$$arg\left(1-\frac{|E|^2i}{3(1-|E|^2)}\right)<arg(1+E^2/3+E^4/5+\ldots)<arg\left(1+\frac{|E|^2i}{3(1-|E|^2)}\right)$$
which 
gives us
$$\arg{\left(1-\frac{3i}{16}\right)}<arg(1+E^2/3+E^4/5+\ldots)<arg{\left(1+\frac{3i}{16}\right)}$$
so certainly
$$-\frac{\pi}{15}<arg(1+E^2/3+E^4/5+\ldots)<\frac{\pi}{15}.$$
Thus 
$$-\frac{\pi}{15}<\arg(\log(\zeta(b)))<\frac{7\pi}{24}+\frac{\pi}{15}=\frac{43\pi}{120}<\frac{3\pi}{8}$$
which puts $\log(\zeta)$ outside the region permitted by the Yoccoz inequality.\\

It remains to consider $b=it$ for $1< t <\infty$. 
First we note that at  $t=1$ (this corresponds to $a=4+3i$) we have 
$|E|=0.563171$ and $\arg(E)=0.0749062$, so $$|\zeta|=|(1+E)/(1-E)|=3.54691.$$
It will suffice to show that for $t\in (1,\infty)$ the value of $|E|\in [0.56,1.0]$ and 
$\arg(E) \in [-0.075,+0.075]$, since $|\zeta|=|(1+E)/(1-E)|$ will then
be greater than $(3+\sqrt{5})/2=2.618\ldots$, the maximum permitted by our Yoccoz inequality.\\

Squaring $E$ gives a formula for $E^2$ as a rational function of $t$:
$$E^2=\frac{b^3(b-4)}{(b^2-2b-2)^2}=\frac{-it^3(it-4)}{(-t^2-2it-2)^2}
=\frac{(t^8+16t^6+36t^4)-i(8t^5-16t^3)}{t^8+16t^6+72t^4+64t^2+16}.$$
It is now easily shown that $|E|$ increases monotonically for $t \in [1,\infty)$, tending to the value of $1$ as $t\to \infty$,
and that $|\arg(E)|$ has its maximum value on $t\in [1,\infty)$ at $t=1$ (in fact, as $t$ 
increases from $1$, 
$\arg(E)$ decreases to a minimum value of about $-0.0388$ and then increases, tending to $0$ as $t\to \infty$:
it passes through $0$ at $t=\sqrt{2}$, which corresponds to $a=3+2\sqrt{2}i$). This completes the proof of Theorem \ref{param-lune}.

\begin{remark}\label{bound}
The bound $(3+\sqrt{5})/2$ on the modulus of the derivative at the fixed point  (established in Corollary \ref{abs_bound} in 
Section \ref{Yoccoz_proof}) is sharp. At the parameter 
value $a=4$ the correspondence ${\mathcal F}_a$ 
is a mating of $z\to z^2-2$ with $PSL(2,{\mathbb Z})$ and so the limit set is connected, indeed it is an interval.
When $a=4$, the fixed point $z_0$ is $-1/\sqrt{5}$ and our formula for
the derivative $\zeta$ of ${\mathcal F}_4$ at $z_0$ gives $\zeta=-(3+\sqrt{5})/2$.  We remark that the classical Yoccoz inequality gives a bound of $2$ for the modulus of the derivative of $ z \to z^2+c$ at the $\alpha$-fixed-point, and that this bound is achieved for $ z \to z^2-2$ at its fixed point $z=-1$.
\end{remark}

{\bf Proof of Corollary \ref{near_7}.} As in our proof of Theorem \ref{param-lune} above, we work in terms of the parameter $b=(a-7)/(a-1)$ in place of $a$. 
In the proof above we established a formula for the multiplier at the alpha-fixed-point, namely:
$$\zeta(b)= \frac{1+E}{1-E}$$
where
$$E=\frac{b(4-b)}{2+2b-b^2}\sqrt{\frac{b}{b-4}},$$
and so for $b$ close to $0$, we have that $E$ is close to zero and 
$$\zeta(b)\sim 1+2E\sim  1+i b^{3/2},$$
where we have neglected all powers of $b$ other than the leading term. Thus for values of $b$ close to the origin, the transformation from the 
$b$-plane to the $(\log{\zeta})$-plane is given by:
$$\log{\zeta}\sim i b^{3/2}.$$
and in the opposite direction it is given by 
$$b \sim \left(i \log{\zeta}\right)^{2/3},$$
that is,
$$a \sim 7 + 6\left( i \log{\zeta}\right)^{2/3},$$ 
and the statement of the Corollary follows.
\qed

\bigskip
\noindent
School of Mathematical Sciences,
Queen Mary University of London, 
London E1 4NS, UK 
\hfill s.r.bullett@qmul.ac.uk

\medskip

\noindent
Instituto de Matem\'atica Pura e Aplicada, 
Estrada Dona Castorina 110, Jardim Bot\^anico, Rio de Janeiro, RJ, CEP 22460-320, Brasil
\hfill luna@impa.br


\end{document}